\documentclass[12pt,draftcls,onecolumn]{IEEEtran}
\IEEEoverridecommandlockouts                              
\usepackage{graphicx}
\usepackage{subfigure}
\usepackage{amsfonts}
\usepackage{amsmath}
\usepackage{amssymb}
\usepackage{psfig}

\newcommand \nn \nonumber
\newcommand \nnd \noindent

\newtheorem{thm}{Theorem}

\newtheorem{lem}{Lemma}

\title{Generalized Voronoi Partition Based Multi-Agent Search using Heterogeneous Sensors}

\author{K.R. Guruprasad
\thanks{K.R. Guruprasad is a Senior Lecturer at the Department of Mechanical Engineering,
National Institute of Technology Karnataka, Surathkal, 575025, India.
{\tt\small krgprao@gmail.com} ({\bf Corresponding author})} and
Debasish Ghose
\thanks{Debasish Ghose
is a Professor at the Guidance, Control, and Decision Systems
Laboratory, Department of Aerospace Engineering, Indian Institute of
Science, Bangalore, 560012, India. {\tt\small
dghose@aero.iisc.ernet.in}}}

\begin{document}

\maketitle 

\begin{abstract} In this paper we propose search strategies for heterogeneous multi-agent systems. Multiple agents, equipped with
communication gadget, computational capability, and sensors having heterogeneous capabilities,
are deployed in the search space to gather information such as presence of targets. Lack of
information about the search space is modeled as an uncertainty
density distribution. The uncertainty is reduced on collection of information
by the search agents. We propose a generalization of  Voronoi
partition incorporating the heterogeneity in sensor capabilities, and  design
optimal deployment strategies for multiple agents, maximizing a
single step search effectiveness. The optimal deployment forms the basis for two search strategies,
namely, {\em heterogeneous sequential deploy and search} and {\em heterogeneous combined deploy
and search}. We prove that the proposed strategies can
reduce the uncertainty density to arbitrarily low level under ideal conditions. We provide a few formal analysis results
related to stability and convergence of the proposed control laws, and to spatial distributedness of the strategies
under constraints such as limit on maximum speed of agents, agents moving with constant speed and limit on sensor
range. Simulation results are
provided to validate the theoretical results presented in the paper.
\end{abstract}


\section{Introduction}
\subsection{Multi-agent systems} Solving complex problems requires higher
intelligence, which nature has gifted to human beings. An
alternative to higher individual intelligence is cooperation among
individuals with limited intelligence. Such group intelligence is
exhibited in nature by swarms of bees, flocks of birds, schools of
fish etc. In these, and in myriad such examples from nature, the
key factor is cooperation with limited, local, and noisy
communication among individuals in a large group. The individuals
are governed by a set of simple behavior leading to a more complex
and useful emergent group behavior. Honey bees' nests, territories
of the male Tilapia Mossambica etc., exhibit a kind of locational
optimization which can be interpreted in terms of centroidal Voronoi
configurations \cite{barlow}. This kind of optimal behavior is also
an outcome of a set of rudimentary decisions by individuals with
local interactions.

Inspired by nature, scientists and engineers have developed the
concept of multi-agent systems with robots, UAVs, etc., as agents.
These multi-agent systems can perform a wide variety of tasks such
as search and rescue, surveillance, achieve and maintain spatial
formations, move as flocks while avoiding obstacles, multiple source
identification and several other tasks.

\subsection{Search using multiple agents} One of the very useful application of multi-agent systems is
search and surveillance. Searching for the presence of targets of interest, survivors in a
disaster, or information of interest in a large, possibly
un-mapped geographical area, is an interesting and practically
useful problem. The problem of searching for targets in unknown
environments has been addressed in the literature in the past
under restrictive conditions \cite{koopman}-\cite{stone}. These
seminal contributions were mostly theoretical in nature and were
applicable to a single agent searching for single or multiple, and
static or moving, targets. Cooperative search by multiple agents have been studied by various
researchers. Enns et al. \cite{enns} use predefined lanes
prioritizing them with the probability of existence of the target.
The vehicles cooperate in that the total path length covered by them
is minimized while exhaustively searching the area. A dynamic
inversion based control law is used to make the vehicles follow the
assigned tracks or lanes while considering the maximum turn radius
constraint. Spires and Goldsmith \cite{spires} use space filling
curves such as Hilbert curves to cover a given space and perform
exhaustive search by multiple robots. Vincent and Rubin
\cite{vincent} address the problem of cooperative search strategies
for unmanned aerial vehicles (UAVs) searching for moving, possibly
evading, targets in a hazardous environment. They use predefined
swarm patterns with an objective of maximizing the target detection
probability in minimum expected time and using minimum number of
UAVs having limited communication range. Beard and McLain
\cite{beard} develop strategies for a team of cooperating UAVs to
visit regions of opportunity without collision while avoiding
hazards in a search area  using dynamic programming methods. The UAVS are also
required to stay within communication range of each other. Flint et
al. \cite{flint} provide a model and algorithm for path planning of
multiple UAVs searching in an uncertain and risky environment, also using
a dynamic programming approach. For this purpose, the search area is divided into
cells and in each cell the probability of existence of a target is
defined. Pfister \cite{pf} uses fuzzy cognitive map to model the
cooperative control process in an autonomous vehicle. In
\cite{poly}-\cite{yang2} the authors use distributed reinforcement
learning and planning for cooperative multi-agent search. The agents
learn about the environment online and store the information in the
form of a {\em search map} and utilize this information to compute
online trajectories. The agents are assumed to be having limitation
on maneuverability, sensor range and fuel. In \cite{yang1} the
authors show a finite lower bound on the search time. Rajnarayan and
Ghose \cite{rajghose} use concepts from team theory to formulate
multi-agent search problems as nonlinear programming problems in a
centralized perfect information case. The problem is then
reformulated in a Linear-Quadratic-Gaussian setting that admits a
decentralized team theoretic solution. Dell et
al. \cite{dell} develop an optimal branch-and-bound procedure with
heuristics such as combinatorial optimization, genetic algorithm and
local start with random restarts, for solving constrained-path
problems with multiple searchers. Sujit and Ghose \cite{sujith} use
concepts of graph theory and game theory to solve the problem of
coordinated multi-agent search. They partition the search space into
hexagonal cells and associate each cell with an uncertainty value
representing lack of information about the cell. As the agents move
through these cells, they acquire information, reducing the
corresponding uncertainty value. Jin et al. \cite{jin} address a
search and destroy mission problem in a military setting with
heterogeneous team of UAVs.

Mobile agents equipped with sensors to gather information about the search area form a sensor network. Optimal deployment of these sensors or agents carrying sensors which is referred to as ``sensor coverage" in the literature, is an important step in achieving effective search. Voronoi partition and its variations are used in sensor network literature. We review the concept of Voronoi partition and some literature on multi-agent systems using this concept.

\subsection{Voronoi partition in sensor network and multi-agent systems}
Voronoi partition (named after Georgy Voronoi \cite{vor2}), also
called {\em Dirichlet tessellation} (named after Gustav Lejeune
Dirichlet \cite{vor1}),  is a widely used scheme of partitioning a
given space based on the concept of ``nearness" of objects such as points in a set
to some finite number of pre-defined locations in the set. In its standard setting Euclidean distance is used as a measure of ``nearness" (see \cite{franz} for a survey). This
concept finds application in many fields such as CAD, image
processing \cite{lvq, vor_gen1} and sensor coverage \cite{bullo1,
bullo2}.
%
There are various generalization of the {\em Voronoi decomposition}
such as {\em weighted Voronoi partitions} and Voronoi partition
based on non-Euclidean metric. The dual of Voronoi diagram is the
{\em Delaunay graph} (named after Boris Delaunay \cite{delaunay}).
These two concepts are very useful in multi-agent search.

A class of problems known as locational optimization (or facility
location) \cite{drezner,okabe}, is used in many applications. These
concepts have been used in sensor coverage literature for optimal
deployment of sensors. Centroidal Voronoi configuration is a standard
solution for this class of problems \cite{du}, where the optimal
configuration of agents is the centroids of the corresponding
Voronoi cells. Cortes et al. \cite{bullo1,bullo2} use these
concepts to solve a spatially distributed optimal deployment
problem for multi-agent systems. A density distribution, as a
measure of the probability of occurrence of an event
along with a function of the Euclidean
distance providing a quantitative assessment of the sensing performance, is used to formulate the problem. Centroidal Voronoi
configuration, with centroid of a Voronoi cell, computed based on
the density distribution within the cell, is shown to be the
optimal deployment of sensors minimizing the sensory error. The
Voronoi partition becomes the natural optimal partitioning due to
monotonic variation of sensor effectiveness with the Euclidean
distance. Schwager et al. \cite{schwager} interpret the density
distribution of \cite{bullo2} in a non-probabilistic framework and
approximate it by sensor measurements. Pimenta et al.
\cite{vkumar} follow a similar approach to address problem with
heterogeneous robots. They let the sensors to be of
different types (in the sense of having different footprints) and
relax the point robots assumptions. Generalization of Voronoi
partition such as power diagrams (or Voronoi diagram in Laguerre
geometry) are used to account for different footprints of the
sensors (assumed to be discs). Due to assumption of finite size of
robots, the robots are assumed to be discs and a free Voronoi
region is defined. A constrained locational optimization problem
is solved.  They also extend the results to non-convex
environments. Ma et al. \cite{ming} use an adaptive triangulation
(ATRI) algorithm based on the Delaunay triangulation
\cite{delaunay}, which is a dual of the Voronoi partition, with
length of the Delaunay edge as a parameter, to achieve non-uniform
coverage.

In \cite{acods}-\cite{krg_thesis} authors address a problem of searching an unknown area with {\em a priori} known uncertainty distribution using multiple agents. The concept of Voronoi partition was used in formulating optimal deployment strategies for multiple agents maximizing the search effectiveness.


\subsection{Motivation and contribution of the paper}
In the literature on multi-agent search, it is largely assumed that the
agents and the sensors are homogeneous in nature. But this
assumption may not be valid in many practical applications. It is
most likely that the sensors will have varied capabilities in
terms of the strength and range, making the problem heterogenous
in nature. We address this issue in this paper, and formulate
and solve a heterogeneous multi-agent search problem. In order to
solve this problem, we present a generalization of the standard
Voronoi partition and use it to design an optimal deployment of
heterogeneous agents.

The generalization of the Voronoi partition proposed in this work takes
into account the heterogeneity in the sensors' capabilities, in
order to design an optimal deployment strategy for heterogeneous
agents. The agent locations are used as sites or
nodes and a concept of a node function, which is the sensor effectiveness
function associated with each node is introduced in place of the usual distance measure. The standard Voronoi partition
and many of its variations can be obtained from this generalization.

In \cite{bullo1,bullo2,schwager} authors use Voronoi partitions for optimal deployment of homogeneous sensors and \cite{vkumar} use power diagrams for the case of heterogeneous agents. In \cite{acods}-\cite{krg_thesis}, authors use Voronoi partition to design multi-agent search strategies for agents with homogeneous sensors. In this paper, we generalize these concepts and incorporate the sensors with heterogeneous capabilities. The optimal deployment
strategy is developed based on the generalized Voronoi partition maximizing the search effectiveness in a given step and forms the basis for two heterogeneous multi-agent search strategies namely, {\em heterogeneous sequential deploy and search} and
{\em heterogeneous combined deploy and search}. We provide convergence results for
the search strategies and also analyze the strategies for spatial
distributedness property. Some preliminary results using heterogeneous
sensors have been earlier reported in \cite{aamas08}. The concepts developed in this work are based on and generalization of those provided in \cite{acods},\cite{isvd07}.

\subsection{Organization of the paper}
The paper is organized as follows. We preview a few mathematical concepts used in this work in Section II. In Section III we
provide a generalization of Voronoi partition. In Section IV we formulate a heterogeneous multi-agent
search problem. The multi-center objective function,
its critical points, the control law responsible for motion of
agents, and its convergence and spatial distributedness property are
also discussed here. In Section V we impose a few constraints on the agents' speeds and provide the convergence proof for the agents' trajectories with the corresponding control laws. In Section VI we propose and analytically study  the {\em heterogeneous sequential
deploy and search} strategy. In Section VII we propose and analyze the {\em heterogeneous combined deploy and search} strategy. We study the effect of limit on sensor range in Section VIII and discuss a few implementation issues in Section IX. Simulation results and discussions are provided in Section X, and finally the paper concludes in Section XI with possible directions for future work.
\section{Mathematical preliminaries} In this section we preview
mathematical concepts such as LaSalle's invariance principle and
Liebniz theorem used in the present work.

\subsection{LaSalle's invariance principle} Here we state LaSalle's
invariance principle \cite{lasalle,lasalle2} used widely to study
the stability of nonlinear dynamical systems. We state the theorem
as in \cite{marquez} (Theorem 3.8 in \cite{marquez}).

Consider a dynamical system in a domain $D$
\begin{equation}
\label{dyn_lasalle} \dot x = f(x)\text{, } f:D \rightarrow
\mathbb{R}^d
\end{equation}

Let $V:D \rightarrow \mathbb{R}$ be a continuously differentiable
function and assume that (i) $M \subset D$ is a compact set, invariant with respect to
the solutions of (\ref{dyn_lasalle}); (ii) $\dot{V} \leq 0$ in $M$; (iii) $E:\{x:x\in M\text{, and } \dot{V}(x) = 0 \}$; that is, $E$
is set of all points of $M$ such that $\dot{V}(x)=0$; and (iv) $N$ is the largest invariant set in $E$.
Then every solution of (\ref{dyn_lasalle}) starting in $M$
approaches $N$ as $ t \rightarrow \infty$.

Here by \emph{invariant set} we mean that if the trajectory is
within the set at some time, then it remains within the set for all
time. Important differences of the LaSalle's invariance principle as
compared to the Lyapunov Theorem are (i) $\dot{V}$ is required to be
negative semi-definite rather than negative definite and (ii) the
function $V$ need not be positive definite (see Remark on Theorem
3.8 in \cite{marquez}, pp 90-91).

\subsection{Leibniz theorem and its generalization}
The Leibniz theorem is widely used in fluid mechanics \cite{kundu},
and shows how to differentiate an integral whose integrand as well
as the limits of integration are functions of the variable with
respect to which differentiation is done. The theorem gives the
formula
\begin{equation}
\label{leibniz} \frac{d}{dy}\int_{a(y)}^{b(y)}F(x,y)dx = \int_a^b
\frac{\partial F}{\partial y}dx + \frac{db}{dy}F(b,y) -
\frac{da}{dy}F(a,y)
\end{equation}

Eqn. (\ref{leibniz}) can be generalized for a $d$-dimensional
Euclidean space as
\begin{equation}
\label{leibnizgen} \frac{d}{dy}\int_{\mathcal{V}(y)}F(x,y)d\mathcal{V}= \int_\mathcal{V}
\frac{\partial F}{\partial y}d\mathcal{V} + \int_{\mathbf{S}}
\mathbf{n}(x).\mathbf{u}(x)FdS
\end{equation}
where, $\mathcal{V} \subset \mathbb{R}^d$ is  the volume in which the
integration is carried out, $d\mathcal{V}$ is the differential
volume element, $\mathbf{S}$ is the bounding hypersurface of $V$,
$\mathbf{n}(x)$ is the unit outward normal to $\mathbf{S}$ and
$\mathbf{u}(x) = \frac{d\mathbf{S}}{dy}(x)$ is the rate at which the surface
moves with respect to $y$ at $x \in \mathbf{S}$.

\section{Generalization of the Voronoi partition} Here we present a generalization
of the Voronoi partition considering the heterogeneity in the
sensors' capabilities. Voronoi partition \cite{vor2,vor1}
is a widely used scheme of partitioning a given space and finds
applications in many fields such as CAD, image processing and sensor
coverage. We can find several extensions or
generalizations of Voronoi partition to suit specific applications
\cite{franz,vor_gen1,okabe}. Herbert and Seidel \cite{herbert}  have introduced an
approach in which, instead of the site set, a finite set of
real-valued functions $f_i: D \mapsto \mathbb{R}$ is used to
partition the domain $D$. Standard Voronoi partition and other known
generalizations can be extracted from this abstract general form.

In this paper we define a generalization of the Voronoi partition
to suit our application, namely the heterogeneous multi-agent
search. We use, (i) the search space as the space to be partitioned, (ii) the site set as the set of points in the search space
which are the positions of the agents in it, and (iii) a set of node functions in place of a distance measure.

Consider a space $Q \subset \mathbb{R}^d$, a set of points called
{\em nodes} or {\em generators} $\mathcal{P} = \{p_1,p_2, \ldots,
p_N \}$, $p_i \in Q$, with $p_i \neq p_j$, whenever $i \neq j$,
and monotonically decreasing analytic functions \cite{krantz} $f_i :
\mathbb{R}^+ \mapsto \mathbb{R}$, where $f_i$ is called a {\em
node function} for the $i$-th node. Define a collection $\{V_i\}$,
$i\in\{1,2,\ldots,N\}$, with mutually disjoint interiors, such that
$Q=\cup_{i \in\{1,2, \ldots, N\}} V_i$, where $V_i$ is defined as
\begin{equation}
\label{vor_fun} V_i = \{ q\in Q  |  f_i(\| p_i - q \|) \geq f_j(\|
p_j - q \|) \quad \forall j \neq i\text{,} \quad j \in \{1,2, \ldots,
N\} \}
\end{equation}

We call $\{V_i\}$, $i\in\{1,2, \ldots, N\}$, as a {\em generalized
Voronoi partition} of $Q$ with nodes $\mathcal{P}$ and node
functions $f_i$. In the standard definition of the Voronoi partition, $f_i(\| p_i - q \|) \geq f_j(\|
p_j - q \|)$ is replaced by $(\| p_i - q \|) \leq (\|
p_j - q \|)$.
Note that,
\begin{itemize}
\item [i)] $V_i$ can be topologically non-connected and may contain other Voronoi cells.

\item [ii)] In the context of multi-agent search problem discussed in this paper, $q
\in V_i$ means that the $i$-th agent is the most effective in
performing search task at point $q$. This is reflected in the $\geq$
sign in the definition. In standard Voronoi partition used for the
homogeneous multi-agent search, $\leq$ sign for
distances ensured that $i$-th agent is most effective in $V_i$
\item  [iii)] The condition that $f_i$ are analytic implies that for every
$i,j \in \{1,2, \ldots, N\}$, $f_i - f_j$ is analytic. By the
property of real analytic functions \cite{krantz},
the set of intersection points between any two node functions is a
set of measure zero. This ensures that the intersection of any two
cells is a set of measure zero, that is, the boundary of a cell is
made up of the union of at most $d-1$ dimensional subsets of
$\mathbb{R}^d$. Otherwise the requirement that the cells should
have mutually disjoint interiors may be violated. Analyticity of
the node functions $f_i$ is a sufficient condition to discount
this possibility.
\item [iv)] The standard Voronoi partition and its generalizations such as multiplicatively and additively weighted Voronoi partitions can be extracted as special cases of the proposed generalization.
\end{itemize}

\begin{thm}
\label{cont_vor} The generalized Voronoi partition depends at least continuously on $\mathcal{P}$.
\end{thm}

\nnd {\it Proof:~} If $V_i$ and $V_j$ are adjacent cells, then all the points $q\in Q$ on the boundary common to them are given by  $\{q \in Q | f_i(\|p_i-q\|)=f_j(\|p_j-q\|)\}$, that is, the intersection of corresponding node functions. Let the $j$-th agent moves by a small distance $dp$. This makes the common boundary between $V_i$ and $V_j$ move by a distance, say $dx$. Now as the node functions are monotonically decreasing and are continuous, it is easy to see that $dx \rightarrow 0$ as $dp \rightarrow 0$. This is true for any pair $i$ and $j$. Thus, the Voronoi partition depends continuously on $\mathcal{P} = \{p_1,p_2,\ldots,p_N\}$. \hfill $\Box$

\subsection*{Generalized Delaunay graph} Delaunay graph is the dual of
Voronoi partition. Two nodes are said to be neighbors (connected by
an edge), if the corresponding Voronoi cells are adjacent. This
concept can be extended to generalized Voronoi partitioning scheme.
For the sake of simplicity we call such a graph a {\em Delaunay
graph}, $\mathcal{G}_D$. Note that the generalized Delaunay graph,
in general, need not have the property of Delaunay triangulation, in
fact, it need not even be a triangulation.

Two nodes are said to be neighbors in a {\em generalized Delaunay
graph}, if the corresponding {\em generalized Voronoi} cells are
adjacent, that is, $(i,j) \in \mathcal{E}_{\mathcal{G}_D}$, the edge
set corresponding to the graph $\mathcal{G}_D$, if $V_i\cap V_j \neq
\emptyset$.

\section{Heterogeneous Multi-Agent Search} In this section we discuss the problem
addressed in this paper. $N$ agents are deployed in the search space $Q \subset \mathbb{R}^d$, a convex polytope, where, lack of information is modeled as an uncertainty density distribution $\phi:Q\mapsto[0,1]$, a continuous function in $Q$. $P(t) = \{p_1(t),p_2(t),\ldots,p_N(t)\}$, $p_i(t) \in Q$ denotes the configuration of the multi-agent system at time $t$, $p_i(t)$ denotes the position
of the $i$-th agent at time $t$. In future, for convenience, we drop
the variable $t$ and refer to the positions by just $p_i$. The agents are assumed to have sensors with varied strength and range, whose effectiveness at a point reduces with distance. The agents get deployed in $Q$, perform search, thereby reducing the uncertainty, and we are looking for optimal utilization of the agents to reduce the uncertainty $\phi(q)$ at each point $q \in Q$ below a specified level.

At each iteration, after deploying themselves {\em optimally}, the
sensors gather information about $Q$, reducing the uncertainty density as,
\begin{equation}
\label{phiupdate_het} \phi_{n+1}(q) = \phi_n(q)\min_i\{\beta_i(\| p_i - q
\|)\}
\end{equation}
where, $\phi_n(q)$ is the density function at the $n$-th iteration;
$\beta: \mathbb{R}\mapsto [0,1]$ is a function of the Euclidian distance
of a given point in space from the agent, and acts as the factor of
reduction in uncertainty by the sensors; and $p_i$ is the position
of the $i$-th sensor. At a given $q \in Q$, only the agent with the
smallest $\beta_i(\| p_i - q \|)$, that is, the agent which can reduce
the uncertainty by the largest amount, is active. If any agent searches within its Voronoi cell, then the updating function (\ref{phiupdate_het}) gets implemented automatically, That is, the function $\min_i\{\beta_i(\|p_i - q \|)\}$ is simply $\beta_i(\|p_i - q \|)$, where $p_i \in V_i$.
Usually, sensors' effectiveness decreases with Euclidean distance,
thus $\beta_i$, which represents the search effectiveness of $i$-th sensor, can be assumed to be a monotonically increasing
analytic function of the Euclidean distance. Equation
(\ref{phiupdate_het}) selects the agent $i$, which is most effective
in performing search task at point $q \in Q$. We will discuss more about the function $\beta_i$ in a later section.

Note that the condition that $\phi_0$ is continuous ensures that $\phi_n$ is continuous for all $n$ as $\beta_i$ are continuous and at any point $q\in Q$, and on the boundary of generalized Voronoi cells corresponding to any pair of agents $i$ and $j$, $\beta_i(\|p_i - q\|) = \beta_j(\|p_j-q\|)$, by the definition of the generalized Voronoi partition.

%
\subsection{Objective function}
Now suppose that the agents have to be deployed in $Q$ in such a way
as to maximize one-step uncertainty reduction, that is, maximize the
effectiveness of one-step of multi-agent search. Consider the objective
function for the $n$-th search step,
\begin{eqnarray}
\label{obj_het_gen}
\mathcal{H}^n &=& \int_Q \Delta\phi_n(q)dQ =  \int_Q (\phi_n(q) - \min_{i}\{\beta_i(\|p_i - q
 \|)\}\phi_n(q))dQ\nn \\
              &=& \sum_{i}\int_{V_{i}}\phi_n(q)(1-\beta_i(\|
p_i - q \|))dQ= \sum_{i}\int_{V_{i}}\phi_n(q)f_i(r_i)dQ
\end{eqnarray}
where, $V_i$ is the {\em generalized Voronoi cell} given by
(\ref{vor_fun}) corresponding to the $i$-th agent, with $f_i(.) =
1-\beta_i(.)$, a monotonically decreasing analytic function as
node function, $\mathcal{P} $ as the set of nodes, and $r_i = \|
p_i - q \|$. Below we provide a result which will be useful in obtaining the critical points of the objective function (\ref{obj_het_gen}).

\begin{lem}
\label{gradient} The gradient of the multi-center objective
function (\ref{obj_het_gen}) with respect to $p_i$ is given by
\begin{equation}
\label{grad} \frac{\partial \mathcal{H}^n}{\partial p_i} =
\int_{V_i}\phi(q)\frac{\partial f_i(r_i)}{\partial p_i}dQ
\end{equation}
where $r_i = \|q-p_i\|$.
\end{lem}

\noindent {\it Proof.~} Let us rewrite (\ref{obj_het_gen}) as
\begin{equation}
\mathcal{H}^n = \sum_{i \in \{1,2,\ldots N\}} \mathcal{H}^n_i
\end{equation}
where $\mathcal{H}^n_i = \int_{V_i}f_i(r_i)\phi(q)dQ$. Now,
\begin{equation}
\frac{\partial \mathcal{H}^n}{\partial p_i} = \sum_{j \in \{1,2,\ldots N\}} \frac{\partial
\mathcal{H}^n_j}{\partial p_i}
\end{equation}

Applying the general form of the Leibniz theorem \cite{kundu}
\begin{eqnarray}
\label{grad_proof} \frac{\partial \mathcal{H}^n}{\partial p_i} &=&
\int_{V_i}\phi(q)\frac{\partial f_i}{\partial p_i}(r_i)dQ  + \sum_{j \in
N_i}\int_{A_{ij}}\mathbf{n}_{ij}(q).\mathbf{u}_{ij}(q)\phi(q)f_i(r_i)dQ\\
\nn
&& + \sum_{j \in
N_i}\int_{A_{ji}}\mathbf{n}_{ji}(q).\mathbf{u}_{ji}(q)\phi(q)f_j(r_j)dQ
\nn
\end{eqnarray}
where, $N_i$ is the set of indices of agents which are neighbors of
the $i$-th  agent in $\mathcal{G}_D$, the generalized Delaunay
graph,  $A_{ij}$ is the part of the bounding surface common to $V_i$
and $V_j$,  $\mathbf{n}_{ij}(q)$ is the unit outward normal to $A_{ij}$ at
$q \in A_{ij}$, $\mathbf{u}_{ij}(q) = \frac{dA_{ij}}{dp_i}(q)$, the rate of
movement of the boundary at $q \in A_{ij}$ with respect to $p_i$.

Note that i) $\mathbf{n}_{ij}(q) = -\mathbf{n}_{ji}(q)$, $\forall q \in A_{ij}$, ii) $\mathbf{u}_{ij}(q) = \mathbf{u}_{ji}(q)$, iii)  $f_i(r_i) = f_j(r_j)$, $\forall q \in A_{ij}$, by definition of the generalized Voronoi partition, and iv) $\phi$
is continuous. Thus, it is clear that for each $j \in N_i$, $
\int_{A_{ij}}\mathbf{n}_{ij}(q).\mathbf{u}_{ij}(q)\phi(q)f_i(r_i)dQ
=
-\int_{A_{ji}}\mathbf{n}_{ji}(q).\mathbf{u}_{ji}(q)\phi(q)f_j(r_j)dQ$, and hence, the last two terms in (\ref{grad_proof}) cancel each other. \hfill $\Box$

\subsection{The critical points}The gradient of the objective
function (\ref{obj_het_gen}) with respect to $p_i$, the location of the
$i$-th node in $Q$, can be determined using (\ref{grad}) (by Lemma
\ref{gradient}) as
\begin{eqnarray}
\label{grad_HLOP} \frac{\partial \mathcal{H}^n}{\partial p_i} &=&
\int_{V_i}\phi(q)\frac{\partial f_i(r_i)}{\partial p_i}dQ =
\int_{V_i}\phi(q)\frac{\partial f_i(r_i)}{\partial
{(r_i)}^2}(p_i-q)dQ\nn \\
&=& -\int_{V_i}\tilde{\phi}(q)(p_i-q)dQ = -\tilde{M}_{V_i}(p_i -
\tilde{C}_{V_i})\nn
\end{eqnarray}
where, $r_i=\|q-p_i\|$ and $\tilde{\phi}(q) =
-\phi(q)\frac{\partial f_i(r_i)}{\partial {(r_i)}^2}$. As $f_i, i \in
\{1,2,\ldots,N\}$ is strictly decreasing,  $\tilde{\phi}(q)$ is
always non-negative. Here $\tilde{M}_{V_i}$ and $\tilde{C}_{V_i}$
are interpreted as the mass and centroid of the cell $V_i$ with
$\tilde{\phi}$ as density. Thus, the critical points are $p_i =
\tilde{C}_{V_i}$, and such a configuration $\mathcal{P}$, of
agents is called a {\em generalized centroidal Voronoi
configuration}.

\begin{thm}
\label{spat_grad_HLOP} The gradient, given by  (\ref{grad_HLOP}),
is spatially distributed over the generalized {\it Delaunay graph}
$\mathcal{G}_D$.
\end{thm}

\noindent {\it Proof.~} The  gradient (\ref{grad_HLOP}) with
respect to $p_i \in \mathcal{P}$, the present configuration,
depends only on the corresponding generalized Voronoi cell $V_i$
and values of $\phi$ and the gradient of $f_i$ within $V_i$. The
Voronoi cell $V_i$ depends only on the neighbors of $p_i$ in $\mathcal{G}_D$. Thus, the
gradient (\ref{grad_HLOP}) can be computed with only local
information, that is, the neighbors of $p_i$ in $\mathcal{G}_D$.
\hfill $\Box$

The critical points are not unique, as with the standard Voronoi
partition. But in the case of a generalized Voronoi partition,
some of the cells could become null and such a condition can lead
to local maxima.

\subsection{Selection of $\beta_i$}
The function $\beta_i: \mathbb{R}\mapsto
[0,1]$ is a {\em sensor detection function} corresponding to the
$i$-th agent. The effectiveness of most sensors decreases with the
Euclidean distance.
Consider
\begin{equation}
\beta_i(r) = 1- k_ie^{-\alpha_i {r_i}^2},\qquad k_i \in (0, 1)\quad
\text{and $\alpha_i > 0$}
\end{equation}
Here, $k_ie^{-\alpha_i {r_i}^2}$ represents the effectiveness of the $i$-th sensor
which is maximum at ${r_i}=0$ and tends to zero as $r \rightarrow
\infty$ and $\beta_i$ is minimum at ${r_i} = 0$ (effecting maximum
reduction in $\phi$) and tends to unity as $r \rightarrow \infty$ (change in $\phi$ reduces to zero as $r$
increases). Most sensors' effectiveness reduces over distance as the signal to noise ratio increases. Thus $\beta_i$, which is upside down Gaussian, can model a wide variety of sensors with two tunable parameters $k_i$ and $\alpha_i$.

%

\subsection{Special cases} Here we discuss a few interesting detection
functions.

\subsubsection*{Case 1: $\beta_i(r_i) = 1 - ke^{-\alpha_i r_{i}^2}$}
\begin{figure}
\centerline{
\subfigure[]{\psfig{figure=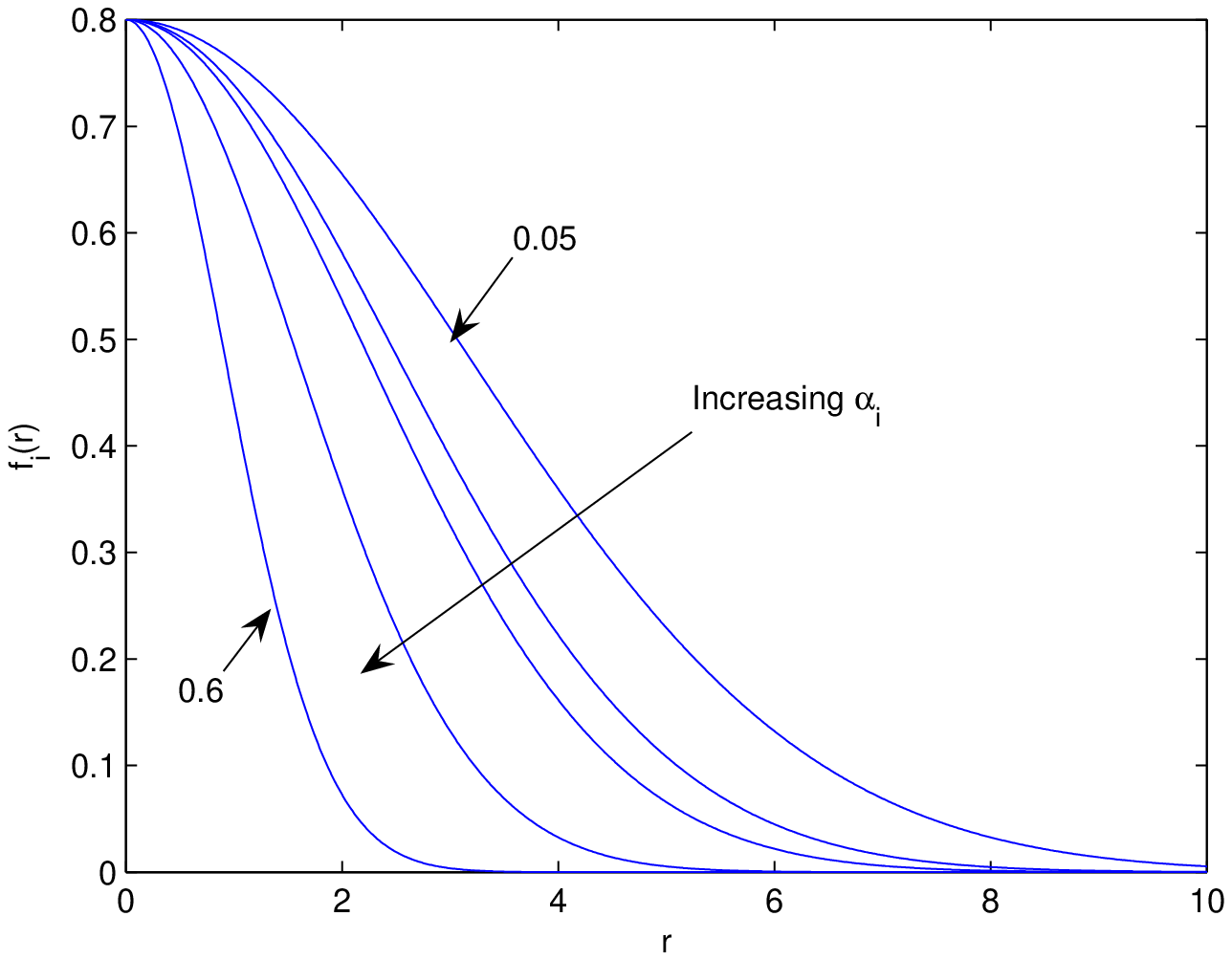,height=6cm,width=6cm}}
\subfigure[]{\psfig{figure=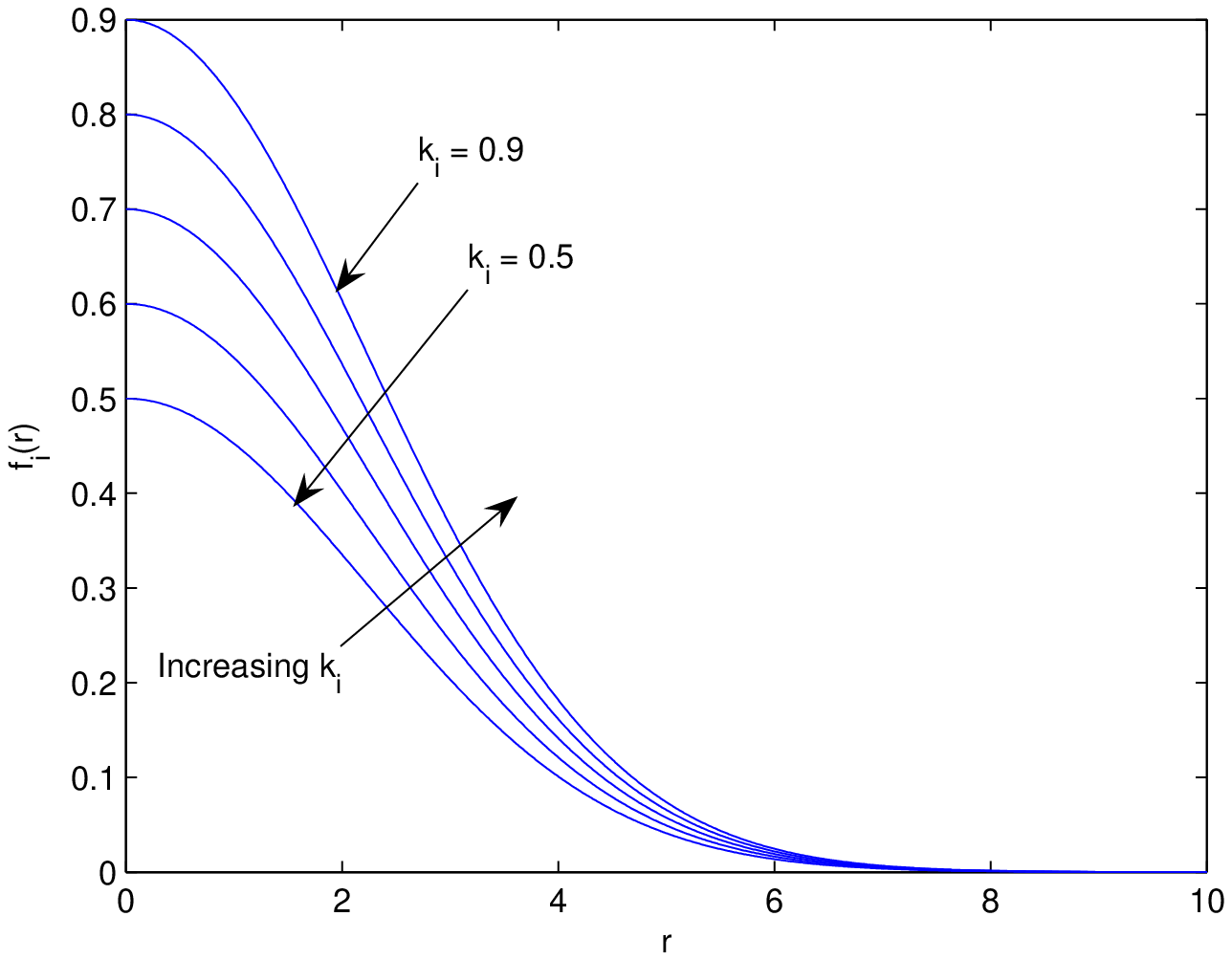,height=6cm,width=6cm}}
\subfigure[]{\psfig{figure=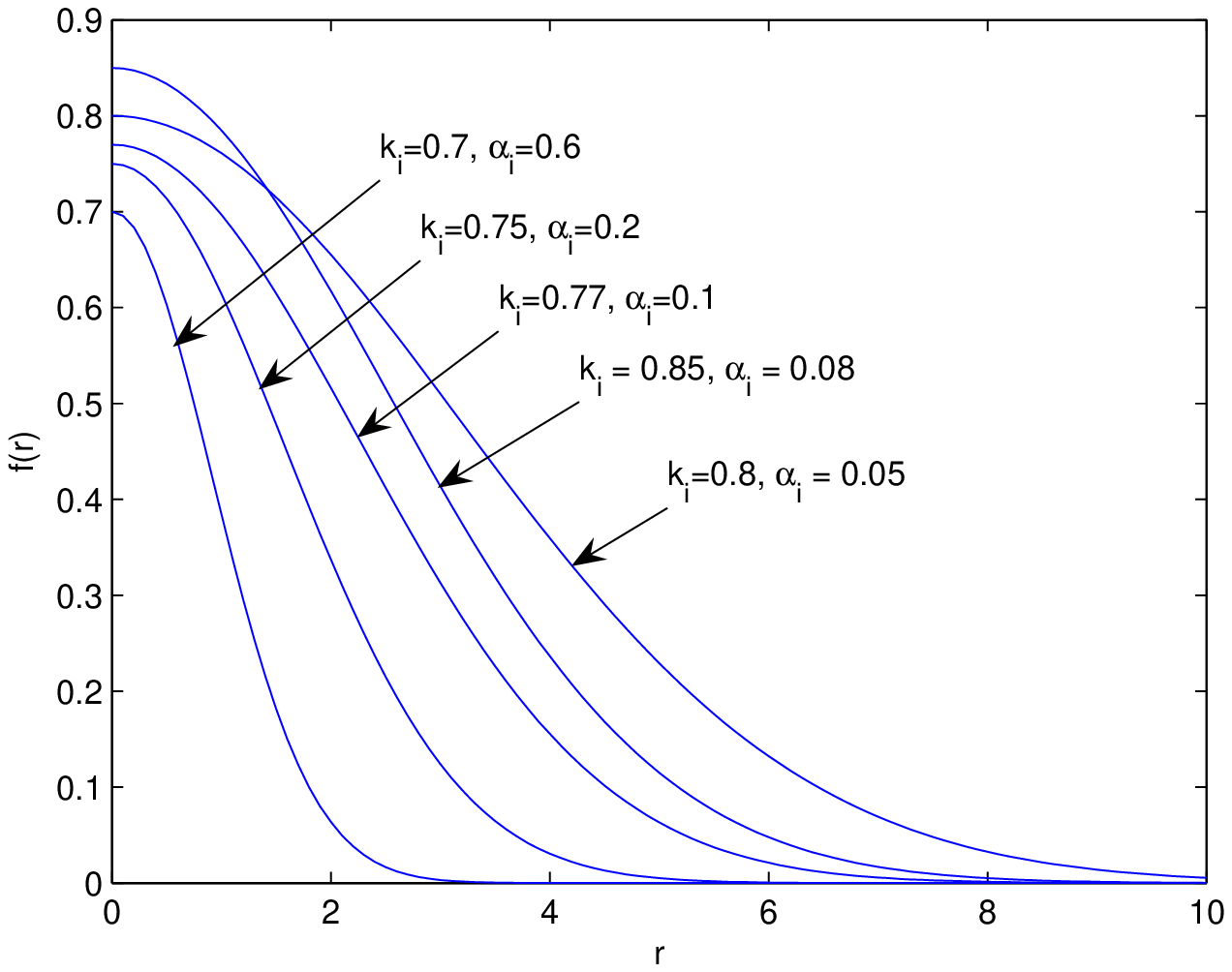,height=6cm,width=6cm}}
}\vspace{-0.08in}
\centerline{
\subfigure[]{\psfig{figure=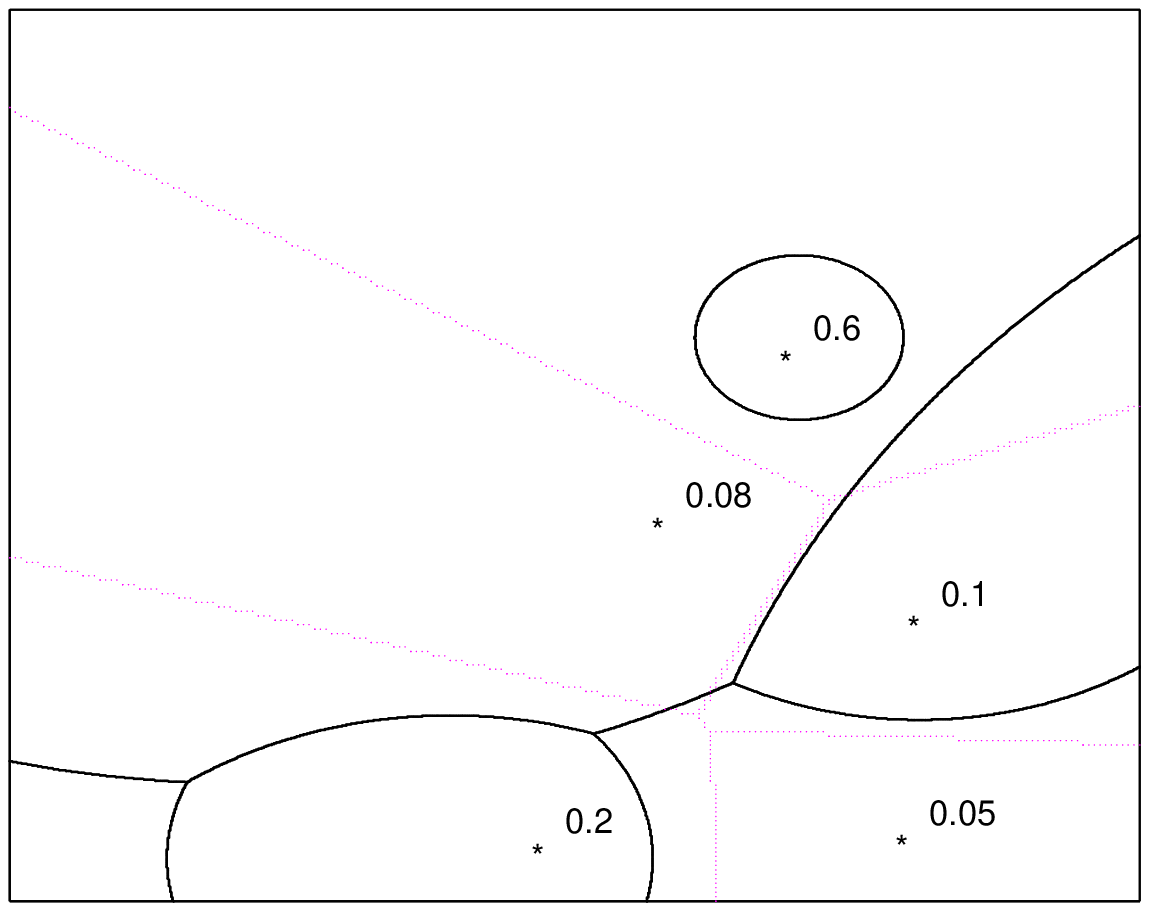,height=6cm,width=6cm}}
\subfigure[]{\psfig{figure=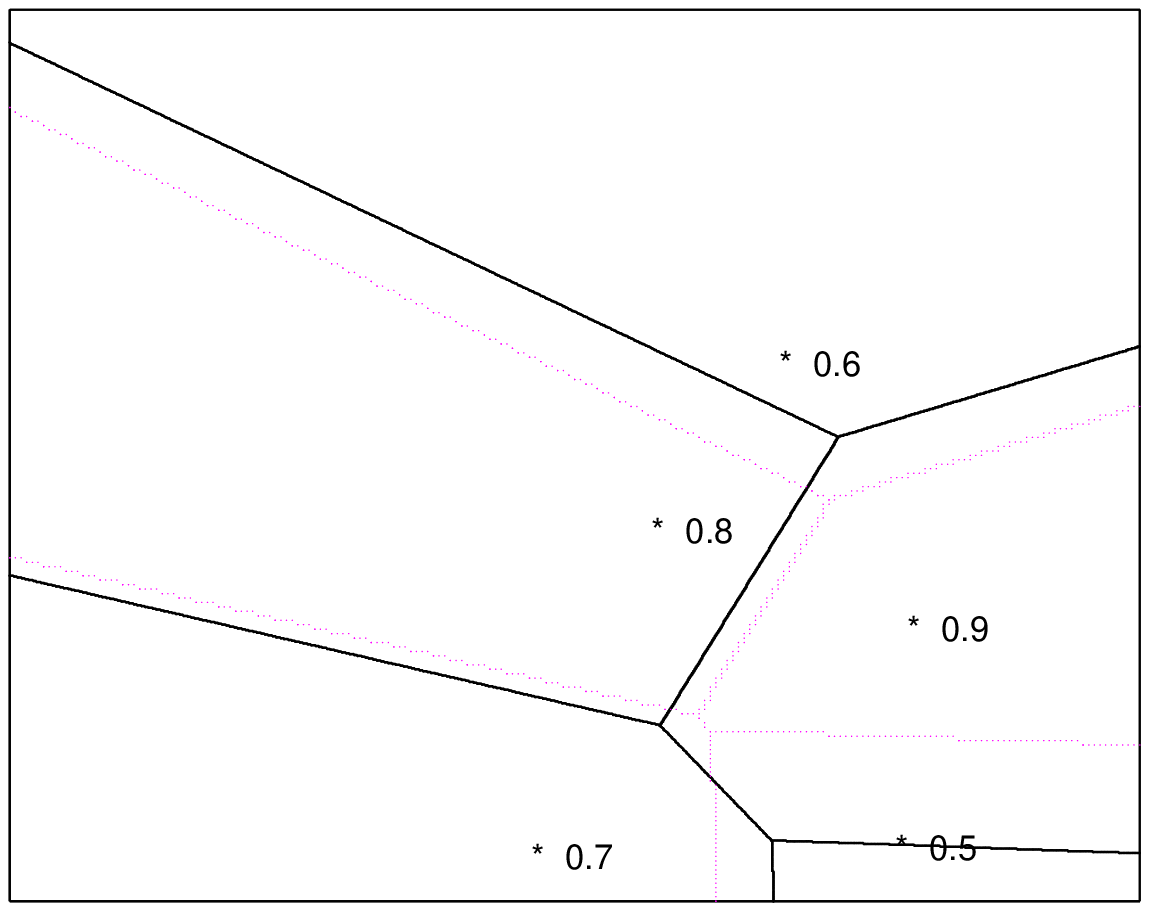,height=6cm,width=6cm}}
\subfigure[]{\psfig{figure=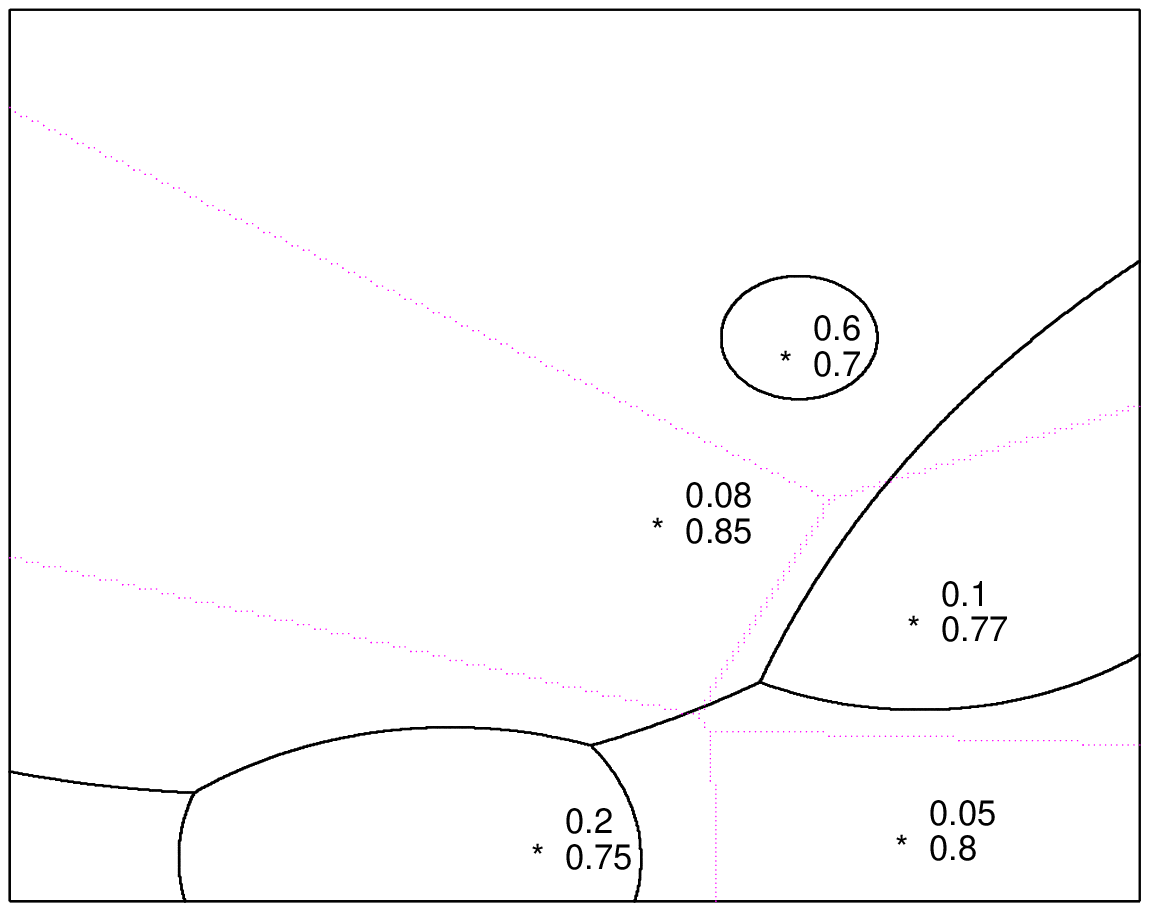,height=6cm,width=6cm}}
}\vspace{-0.2in}\caption{The sensor effectiveness $\beta_i(r) = 1 - k_ie^{-\alpha_i r_{i}^2}$ with a) varying $\alpha_i$, b) varying $k_i$, and c) both $\alpha_i$ and $k_i$ varying. d), e), and f) are the corresponding Voronoi partitions with
'*' showing the locations of nodes or sites, the numbers indicating the parameters, and dotted lines showing the standard Voronoi partition.}\label{betacases}

\end{figure}

\noindent Here we consider exponential sensor effectiveness and
assume that the parameter $\alpha$ is different for different
sensors while $k$ remains the same for all of them, that is, all
the agents have sensors with the same maximum effectiveness with
different sensor reach (Figure \ref{betacases} (a)). This case leads to
multiplicatively weighted Voronoi partition. The Voronoi cells
can be non-connected and also can have one or more Voronoi cells
embedded within a cell. Within $V_i$, $\tilde{\phi}_n(\cdot) =
-\alpha_i\phi_n(\cdot)ke^{-\alpha_i r_{i}^2}$.

Figure \ref{betacases} (d) shows a Voronoi partition for this case with
$k=0.8$. The parameter $k$ does not affect the Voronoi partition. It
is easy to show that the partition is a multiplicatively weighted
Voronoi partition with weights $\sqrt{\alpha_i}$. The Voronoi cell
corresponding to node with $\alpha_i = 0.6$, the highest $\alpha_i$
among all the nodes, implying that it is the least effective, is
embedded within the cell corresponding to the node with $\alpha_i =
0.08$, the smallest value, and hence the most effective. The Voronoi
cell corresponding to the node with $\alpha_i =0.05$ is actually two
cells separated by the cell corresponding to node with $\alpha_i =
0.2$.

\subsubsection*{Case 2. $\beta_i(r_i) = 1 - k_ie^{-\alpha r_{i}^2}$}
Here we let the parameter $k_i$ be different for different sensors
while keeping the value of $\alpha$ same. The agents in this case
have sensors with varying maximum effectiveness (Figure \ref{betacases} (b)). Within $V_i$, $\tilde{\phi}_n(\cdot) =
-\alpha\phi_(\cdot)k_ie^{-\alpha r_{i}^2}$.

\begin{thm}
For $Q \subset \mathbb{R}^2$, the boundaries of Voronoi cells
corresponding to the case $\beta_i(r) = 1 - k_ie^{-\alpha r^2}$
are straight line segments.
\end{thm}

\nnd {\it Proof.~} Let $q = (x,y)'$ be a point on the boundary of
cells $V_i$ and $V_j$. Then we have,
\begin{equation}
\label{case2} k_1e^{-\alpha (\parallel p_1 - q \parallel)^2} =
k_2e^{-\alpha (\parallel p_2 - q \parallel)^2}
\end{equation}

Now let $p_1 = (p_{11}, p_{12})'$ and $p_2 = (p_{21},p_{22})'$.
Solving for $q$ from (\ref{case2}) gives,
\begin{equation}
\label{line} y = \frac{(p_{11} - p_{21})}{(p_{22} - p_{12})}x +
\frac{(1/\alpha)ln(k_1/k_2) + (p_{21}^2 + p_{22}^2) -
(p_{11}^2 + p_{12}^2)}{2(p_{22} - p_{12})}= M x + C
\end{equation}
where, $M$ and $C$ are constants for a given $\mathcal{P}$.
Eqn. (\ref{line}) is true for any $q$ on the boundary of cells
$V_i$ and $V_j$ and hence, the boundary is a straight line
segment. \hfill$\Box$
%

Figure \ref{betacases}(e) shows the corresponding Voronoi partition
along with the standard Voronoi partition. Each cell is made up of
line segments.

\subsubsection*{Case 3.  $\beta_i(r_i) = 1 - k_ie^{-\alpha_i r_{i}^2}$}
Here we vary both the parameters. This case is in the general
category of Voronoi partitions (\ref{vor_fun}). Within $V_i$, $\tilde{\phi}_n(\cdot) =
-\alpha_i\phi_n(\cdot)k_ie^{-\alpha_i r_{i}^2}$
%

Figure \ref{betacases} (f) shows the corresponding Voronoi partition
along with the standard Voronoi partition, while Figure \ref{betacases} (c) shows corresponding sensor effectiveness.

\subsubsection*{Case 4.  $1-\beta_i(r_i) = -r_{i}^2$}
Here we relax the condition that $\beta_i(\cdot) \in [0,1]$. This
will lead to standard centroidal Voronoi configuration
\cite{bullo1,bullo2}.

\subsection{The control law} The critical points of the objective function (\ref{obj_het_gen}) are the respective centroids. Here we discuss a control law making the agents move toward respective centroids.

Typically search problems do not consider dynamics of search agents as the focus is more on the effectiveness
of search, that is, being able to identify region of high uncertainty and distribute search effort to reduce uncertainty.
Moreover, it is usually assumed that the search region is large compared to the physical size of the agent or the
area needed for the agent to maneuver. We assume the first order dynamics for agents to demonstrate the search strategy presented in this work.
Let us consider the system dynamics as
\begin{equation}
\label{dyn1_het} \dot p_{i} = u_i
\end{equation}
Consider the control law
\begin{equation}
\label{ctrl1_het} u_i = -k_{prop}(p_{i} - \tilde{C}_{V_{i}})
\end{equation}
Control law (\ref{ctrl1_het}) makes the agents move toward
$\tilde{C}_{V_{i}}$ for positive gain $k_{prop}$.

It is not necessary that $\tilde{C}_{V_i} \in V_i$, but $\tilde{C}_{V_i} \in Q$ is
true always and this fact ensures that Q is an invariant set for
(\ref{dyn1_het}) under (\ref{ctrl1_het}).

\noindent {\it Remark 2.}~ It can be seen that the control law
(\ref{ctrl1_het}) is spatially distributed over the generalized Delaunay graph
$\mathcal{G}_D$.

\begin{thm}
\label{LaSalle} The trajectories of the agents governed by the
control law (\ref{ctrl1_het}), starting from any initial condition
$\mathcal{P}(0) \in Q^N$, will asymptotically converge to the
critical points of $\mathcal{H}^n$.
\end{thm}

\noindent {\it Proof~}. Here we will use LaSalle's invariance discussed earlier.
Consider $V(\mathcal{P}) = -\mathcal{H}^n$.
\begin{equation}
\label{vdot}
\begin{array}{lcl}
{\dot V}(\mathcal{P}) &=& -\frac{d\mathcal{H}^n}{dt}=
 -\sum_i \frac{\partial \mathcal{H}^n}{\partial p_{i}}\dot{p}_{i}= 2\alpha \sum_i \tilde{M}_{V_{I}}(p_{i} -
 \tilde{C}_{V_{i}})(-k_{prop}(p_{i} -
C'_{V_{i}}))\\
&=& -2\alpha k_{prop}\sum_i  \tilde{M}_{V_{i}}(p_{i} -
\tilde{C}_{V_{i}})^2
\end{array}
\end{equation}

We observe that $V: Q\mapsto \mathbb{R}$ is continuously differentiable in
$Q$ by Theorem \ref{cont_vor}; $M = Q$ is a compact invariant set; ${\dot V}$ is negative definite in $M$; $E = \dot{V}^{-1}(0) = \{\tilde{C}_{V_{i}}\}$ which itself is the largest invariant subset of $E$ by the control law
(\ref{ctrl1_het}). Thus by LaSalle's invariance principle, the trajectories of the
agents governed by control law (\ref{ctrl1_het}), starting from any
initial configuration $\mathcal{P}(0) \in Q^N$, will asymptotically
converge to set $E$, the critical points of $\mathcal{H}^n$. \hfill
$\Box$

It can be noted that the centroid is computed based on the density information. The generalized Voronoi partition is updated
as the agents move and the centroids are recomputed. At the end of a deployment step, the control law (\ref{ctrl1_het}) ensures
that each agent is at (or sufficiently close to) the centroid of the corresponding generalized Voronoi cell, guaranteeing maximal
uncertainty reduction. It is known that the gradient descent/ascent is not guaranteed to find a global optimal solution (see \cite{bullo1} and references therein). Thus, we can only guarantee a local optimum to the optimization problem.

 To implement the control law, centroid of each generalized Voronoi cell needs to be computed.  The computational overhead of computing the centroid can be reduced at the cost of slower convergence using methods reported in the literature such as random sampling and stochastic approximation \cite{lvq,pages}. In addition, we discretize the search space into grids while implementing the strategy. This simplifies the computation of the centroid of generalized Voronoi cells. There are a few efficient algorithms implementing computation related to standard Voronoi partition. Addressing the computational related issues and development of efficient algorithms for computing the generalized Voronoi partitions will help implementating the search strategies presented in this paper more effectively.  The main focus of this work is design and demonstration of heterogenous multi-agent search strategies and finer issues such as computation complexities are beyond the scope of this paper.

\section{Constraints on agents' speed} We proposed a control law to guide the agents toward the critical points, that is, to their respective centroid, and observed that the closed loop system for agents modeled as first order dynamical system, is globally asymptotically stable. Here we impose some of constraints on the agent speeds and analyze the dynamics of closed loop system.

\subsection{Maximum speed constraint} Let the agents have a constraint
on maximum speed of ${U_{max}}_i$, for $i=1,\ldots,n$. Now consider
the control law
\begin{equation}
\label{ctrl_with_sat} u_i =
\begin{cases}
-k_{prop}(p_{i} - \tilde{C}_{V_{i}}) & \text{If $u_i \leq {U_{max}}_i$} \\
-{U_{max}}_i\frac{(p_{i} - \tilde{C}_{V_{i}})}{\parallel(p_{i} -
\tilde{C}_{V_{i}})\parallel} & \text{Otherwise}
\end{cases}
\end{equation}

The control law (\ref{ctrl_with_sat}) makes the agents move toward
their respective centroids with saturation on speed.

\begin{thm}
\label{saturation_speed_stability} The trajectories of the agents
governed by the control law (\ref{ctrl_with_sat}), starting from any
initial condition $\mathcal{P}(0) \in Q^N$, will asymptotically
converge to the critical points of $\mathcal{H}^n$.
\end{thm}

\noindent {\it Proof.~} Consider  $V(\mathcal{P}) = -\mathcal{H}^n$.
\begin{equation}
\label{vdot_sat}
\begin{array}{lcl}
{\dot V}(\mathcal{P}) &=& -\frac{d\mathcal{H}^n}{dt} = -\sum_{i\in\{1,2,\ldots,N\}} \frac{\partial \mathcal{H}^n}{\partial p_{i}}\dot{p}_{i}\\
&=& \begin{cases} 2\alpha \sum_{i\in\{1,2,\ldots,N\}} \tilde{M}_{V_{I}}(p_{i} -
 \tilde{C}_{V_{i}})(-k_{prop})(p_{i} - \tilde{C}_{V_{i}}) \text{, ~~If $u_i \leq {U_{max}}_i$} \\
2\alpha \sum_{i\in\{1,2,\ldots,N\}} \tilde{M}_{V_{I}}(p_{i} -
 \tilde{C}_{V_{i}})(-{U_{max}}_i)\frac{(p_{i} -
C'_{V_{i}})}{(\|p_{i} - \tilde{C}_{V_{i}}\|)}\text{, ~~otherwise}\\
\end{cases} \\
&=&
\begin{cases}
-2\alpha k_{prop}\sum_{i\in\{1,2,\ldots,N\}}  \tilde{M}_{V_{i}}(\|p_{i} -
\tilde{C}_{V_{i}}\|)^2  \text{, ~~If $u_i \leq {U_{max}}_i$}  \\
-2\alpha \sum_{i\in\{1,2,\ldots,N\}} {U_{max}}_i\tilde{M}_{V_{i}}\frac{(\|p_{i} -
\tilde{C}_{V_{i}}\|)^2}{(\parallel p_{i} - \tilde{C}_{V_{i}}\parallel)}
 \text{, ~~ otherwise}
\end{cases}
\end{array}
\end{equation}

We observe that $V: Q\mapsto \mathbb{R}$ is continuously differentiable in
$Q$ as $\{V_i\}$ depends at least continuously on $\mathcal{P}$ (Theorem \ref{cont_vor}), and $\dot V$ is continuous as $u$ is continuous; $M = Q$ is a compact invariant set; ${\dot V}$ is negative definite in $M$; $E = \dot{V}^{-1}(0) = \{\tilde{C}_{V_{i}}\}$, which  itself is the largest invariant subset of $E$ by the control law (\ref{ctrl_with_sat}). Thus, by LaSalle's invariance principle, the trajectories of the
agents governed by control law (\ref{ctrl_with_sat}), starting from
any initial configuration $\mathcal{P}(0) \in Q^N$, will
asymptotically converge to the set $E$, the critical points of
$\mathcal{H}^n$, that is, the generalized centroidal Voronoi configuration with
respect to the density function as perceived by the sensors. \hfill
$\Box$

\subsection{Constant speed control} The agents may have a constraint
of moving with a constant speed $U_i$. But we let the agents slow down as they approach the critical points. For $i=1,\ldots,n$, consider
the control law
\begin{equation}
\label{const_speed_ctrl} u_i =
\begin{cases} -U_i\frac{(p_{i} -
\tilde{C}_{V_{i}})}{\parallel(p_{i} - \tilde{C}_{V_{i}})\parallel}
\text{, if $\|p_i - \tilde{C}_{V_i}\| \geq 1$} \\
-U_i(p_i - \tilde{C}_{V_i}) \text{, otherwise} \end{cases}
\end{equation}

The control law (\ref{const_speed_ctrl}) makes the agents move
toward their respective centroids with a constant speed of $U_i
> 0$ when they are far off from the centroids and slow down as they approach them.

\begin{thm}
\label{const_speed_stability} The trajectories of the agents
governed by the control law (\ref{const_speed_ctrl}), starting from
any initial condition $\mathcal{P}(0) \in Q^N$, will asymptotically
converge to the critical points of $\mathcal{H}^n$.
\end{thm}

\noindent {\it Proof}. Consider  $V(\mathcal{P}) = -\mathcal{H}^n$,
where $\mathcal{P} = \{p_{1}, p_{2}, ... ,p_{N}\}$ represents the
configuration of $N$ agents.
\begin{equation}
\label{vdot_const}
\begin{array}{lcl}
{\dot V}(\mathcal{P}) &=& -\frac{d\mathcal{H}^n}{dt} = -\sum_{i\in\{1,2,\ldots,N\}} \frac{\partial \mathcal{H}^n}{\partial p_{i}}\dot{p}_{i}\\
&=& 2\alpha \sum_{i\in\{1,2,\ldots,N\}} \tilde{M}_{V_{i}}(p_{i} -
 \tilde{C}_{V_{i}})(-U_i)\frac{(p_{i} -
\tilde{C}_{V_{i}}))}{(\parallel p_{i} - \tilde{C}_{V_{i}}\parallel)}\\
&=& \begin{cases} -2\alpha \sum_{i\in\{1,2,\ldots,N\}}U_i
\tilde{M}_{V_{i}}\frac{(\|p_{i} -
\tilde{C}_{V_{i}}\|)^2}{(\parallel p_{i} -
\tilde{C}_{V_{i}}\parallel)} \text{,~~ if $\|p_i - \tilde{C}_{V_i}\| \geq 1$} \\
-2\alpha \sum_{i\in\{1,2,\ldots,N\}}U_i \tilde{M}_{V_{i}}(p_i - \tilde{C}_{V_i}) \text{, ~~otherwise} \end{cases}
\end{array}
\end{equation}

We observe that $V: Q\mapsto \mathbb{R}$ is continuously differentiable in
$Q$ as $\{V_i\}$ depends at least continuously on $\mathcal{P}$ (Theorem \ref{cont_vor}), and $\dot V$ is continuous as $u$ is continuous; $M = Q$ is a compact invariant set; ${\dot V}$ is negative definite in $M$; $E = \dot{V}^{-1}(0) = \{\tilde{C}_{V_{i}}\}$, which itself is the largest invariant subset of $E$ by the control law
(\ref{const_speed_ctrl}). Thus, by LaSalle's invariance principle, the trajectories of the
agents governed by control law (\ref{const_speed_ctrl}), starting
from any initial configuration $\mathcal{P}(0) \in Q^N$, will
asymptotically converge to the set $E$, the critical points of
$\mathcal{H}^n$, that is, the generalized centroidal Voronoi configuration with
respect to the density function as perceived by the sensors. \hfill
$\Box$

\section{Heterogenous sequential deploy and search (HSDS)} In this strategy, the
agents are first deployed optimally according to the objective
function (\ref{obj_het_gen}) and the search task is performed
reducing the uncertainty density at the end of the deployment
step. This iteration of ``deploy" and ``search" in a sequential manner
continues till the uncertainty density is reduced below a required
level. The iteration count $n$ in (\ref{phiupdate_het}) refers to
the number of ``deploy and search" steps. The control law
(\ref{ctrl1_het}) is used to make the agents move toward the
critical points, that is, the centroids of the corresponding
cells.

\noindent {\it Remark 3.~} It is straightforward to prove that
the {\em heterogeneous sequential deploy and search} strategy is
spatially distributed over the generalized Delaunay graph $\mathcal{G}_D$.

\begin{thm}
\label{DNS_convergence_het} The {\em heterogeneous sequential
deploy and search} strategy can reduce the average uncertainty to
any arbitrarily small value in a finite number of iterations.
\end{thm}

\noindent {\it Proof.~} Consider the uncertainty density update
law (\ref{phiupdate_het}) for any $q \in Q$,
\begin{equation}
\label{phiupdate_n-1} \phi_n(q) = (1 - k_ie^{-\alpha_i
{r_i}^2})\phi_{n-1}(q)= \gamma_{n-1}\phi_{n-1}(q)
\end{equation}
where, $r_i$ is the distance of point $q \in Q$ from the $i$-th
agent, such that $q \in V_i$, the Voronoi partition corresponding
to it and define $\gamma_{n-1} = (1 - k_ie^{-\alpha_i {r_i}^2})$. Note that in HSDS, the agents are located in respective centroid at the time of performing search.

Applying the above update rule recursively, we have,
\begin{equation}
\label{phiupdate_recur} \phi_n(q) = \gamma_{n-1}\gamma_{n-2}\ldots
\gamma_1\gamma_0\phi_0(q)
\end{equation}

Let $D(Q) = \max_{p,q \in Q}(\parallel p - q \parallel)$. It
should be noted that
\begin{itemize}
\item[1.] $0 < k_i < 1$
\item[2.] $0 \leq r_i \leq D(Q)$. $D(Q)$ is bounded since the set $Q$
is bounded.
\item[3.] $0 \leq \gamma_j \leq 1 - ke^{-\alpha \{D(Q)^2\}} = l$ (say), $j \in \mathbb{N}$;
and $l < 1$, where $k = \min_i(k_i)$ and $\alpha =
\max_i(\alpha_i)$
\end{itemize}
Now consider the sequence $\{\Gamma\}$ ,
\begin{displaymath}
\Gamma_n = \gamma_n\gamma_{n-1}\ldots \gamma_1\gamma_0  \leq
l^{n+1}
\end{displaymath}
Taking limits as $n \rightarrow \infty$,
\begin{displaymath}
\lim_{n \rightarrow \infty}\Gamma_n \leq \lim_{n \rightarrow
\infty}l^{n+1} = 0
\end{displaymath}
Thus,
\begin{displaymath}
\lim_{n \rightarrow \infty}\phi_{n}(q) = \lim_{n \rightarrow
\infty}\Gamma_{n-1}\phi_0(q) = 0
\end{displaymath}

As the uncertainty density $\phi$ vanishes at each point $q \in Q$
in the limit, the average uncertainty density over $Q$ is bound to
vanish in the limit as $n \rightarrow \infty$. Thus, the HSDS
strategy can reduce the average uncertainty to any arbitrarily
small
value in a finite number of iterations.\hfill $\Box$\\

It should be noted that the above proof does not depend on the control law. The theorem depends only on the outcome of the choice of the updating function (\ref{phiupdate_het}), along with the fact that there is no limit on sensor range  and the search space $Q$ is bounded. In addition, it does not address the issue of the optimality of the strategy which, in fact, depends on the control law responsible for the motion of the agents. However, in HSDS, the reduction in uncertainty in each ``deploy and search" step is maximized. The reduction in the uncertainty in each step in HSDS is
\begin{equation}
\label{HSDSrate}
\mathcal{H}_n^* = \sum_i \int_{V_i} \phi_n(q)k_ie^{-\alpha_i(\|\tilde{C}_{V_i} - q\||)^2} dQ
\end{equation}
which is the maximum possible reduction in a single step. The deployment is such that uncertainty will be reduced
to a maximum possible extent in a step, given by the above formula.

\section{Heterogeneous combined deploy and search (HCDS)} Here
we propose a {\em heterogeneous combined deploy and search} (HCDS)
strategy, where, instead of waiting for the completion of optimal
deployment, as in HSDS, agents perform search as they are moving
toward the respective centroids in discrete time intervals.

\subsection{Density update} Here we provide the problem formulation for
the {\em heterogeneous combined deploy and search} strategy. Assume that the
index $n$ represents the intermediate step at which the search is performed and uncertainty
density is updated. Using the uncertainty density update rule (\ref{phiupdate_het}) discussed earlier we can get,
\begin{equation}
\label{ContSrchDiffrence} \Delta_n \phi(q) = \phi_{n+1}(q) -
\phi_n(q) = \phi_n(q) \min_i (1 - \beta_i(\parallel p_i - q
\parallel))
\end{equation}
Define,
\begin{equation}
\label{bigphi} \Phi_n = \int_Q \phi_n(q)dQ
\end{equation}
Integrating (\ref{ContSrchDiffrence}) over $Q$,
\begin{equation}
\label{contsrch2} \Delta \Phi_n = \sum_{i\in\{1,2,\ldots,N\}}
\int_{V_i}\phi_n(q)(1 - \beta_i(\parallel p_i - q
\parallel))dQ
\end{equation}

\subsection{Objective function}
The objective function (\ref{obj_het_gen}), used for {\em heterogeneous sequential deploy
and search} strategy, is fixed for each deployment step as $\phi_n(q)$ is
fixed for the $n$-th iteration. In {\em heterogeneous combined deploy and search},
the search task is performed as the agents move. The following objective function is maximized in order to maximize the uncertainty reduction at the $n$-th search step,
\begin{equation}
\label{obj_cs1_het} \mathcal{H}^n = \Delta \Phi_n
               = \sum_{i\in\{1,2,\ldots,N\}}
\int_{V_i}\phi_n(q)(1 - \beta_i(\parallel p_i - q
\parallel))dQ
\end{equation}
Note that the above objective function is same the as (\ref{obj_het_gen}) except for the fact that $n$ in this case represents the search step count, whereas in (\ref{obj_het_gen}) it represents ``deploy and search" step count.
For $\beta_i(r) = 1- k_ie^{-\alpha_i r^2}$, the objective function
(\ref{obj_cs1_het}) becomes,
\begin{equation}
\label{obj_cs2_het} \mathcal{H}^n
=\sum_{i\in\{1,2,\ldots,N\}}
\int_{V_i}\phi_n(q)k_ie^{-\alpha_i r_i^2}dQ
\end{equation}

It can be noted that for a given $n$, the uncertainty density $\phi_n(q)$ at any $q \in Q$ is constant. Thus, the gradient of the objective function (\ref{obj_cs2_het}) with respect to $p_i$ can be computed as in HSDS. The gradient is given as,
\begin{eqnarray}
\label{grad_cs_het}
\frac{\partial\mathcal{H}^n}{\partial p_i} &=&
\sum_{i\in\{1,2,\ldots,N\}}\int_{V_{i}}\phi_n(q)k_ie^{-\alpha_i (\| p_i - q\|)^2}(-2\alpha_i)(p_{i}-q)dQ \nn \\
&=& -2\alpha_i \tilde{M}_{V_{i}}(p_{i} - \tilde{C}_{V_{i}})
\end{eqnarray}
As in the case of HSDS, the critical points for any given $n$ are $p_i = \tilde{C}_{V_i}$. But the uncertainty changes in every time step and hence the critical points also change. Hence, the corresponding critical points are only the instantaneous critical points. As the instantaneous critical points  of the objective function (\ref{obj_cs2_het}) are similar to those of (\ref{obj_het_gen}), we  use (\ref{ctrl1_het}), the
control law used for HSDS strategy.

The instantaneous critical points and the gradient (\ref{grad_cs_het}) are used in control law (\ref{ctrl1_het}) only to make the agents move toward the instantaneous centroid rather than deploying them optimally. Thus, it is not possible to prove any optimality of deployment and we do not prove the convergence of the trajectories here in case of HCDS. Agents perform more frequent searches instead of waiting till the optimal deployment maximizing per step uncertainty reduction.

\noindent {\it Remark 4.~} It is straightforward to show that the
continuous time {\em heterogeneous combined deploy and search}
strategy is spatially distributed over the generalized Delaunay graph
$\mathcal{G}_D$.

\begin{thm}
\label{CS_convergence_cont_het} The {\em
heterogeneous combined deploy and search} strategy can reduce the
average uncertainty to any arbitrarily small value in finite time.
\end{thm}

\noindent{\it Proof.~} The proof is similar to Theorem
\ref{DNS_convergence_het} as the density update law is the same. The
differences are that, i) in  {\em heterogeneous combined deploy
and search}, the density update occurs every time step and $n$ represents search step count rather than `deploy and search' count, and ii) the agent configuration in Eqn. (\ref{phiupdate_n-1}) need not be optimal. Even when the agent configuration is non-optimal, $\gamma_n$ are strictly less than unity as noted in the proof of Theorem \ref{DNS_convergence_het}. \hfill $\Box$

As in the case of HSDS, it should be noted that the above proof does not depend on the control law. The theorem depends only on the outcome of the choice of the updating function (\ref{phiupdate_het}), along with the fact that there is no limit on sensor range  and the search space $Q$ is bounded. In addition, it does not address the issue of the optimality of the strategy which, in fact, depends on the control law responsible for the motion of the agents. In HCDS, instead of waiting till the single step uncertainty reduction is maximized, agents perform frequent searches. Though the amount of uncertainty in each step is less than that in HSDS, increased instances of search ensure faster reduction in uncertainty density.

\section{Limited range sensors} We have not considered any limitation on the sensor range in formulating the multi-agent search strategies in previous sections. But in reality the sensors will have
limited range. In this section we formulate search problem for heterogeneous agents having limit on their sensor range.

Let $R_i$ be the limit on range of the sensors and
$\bar{B}(p_i,R_i)$ be a closed ball centered at $p_i$ with a radius
of $R_i$. The $i$-th sensor has access to information only from
points in the set $V_i \cap \bar{B}(p_i,R_i)$. Let us also assume
that $f_i(R_i) = f_j(R_j)\text{, } \forall i\text{,}j \in
\{1,\ldots,N\}$, that is, we assume that the cutoff range for all
sensors is same. Consider the objective function to be maximized,
\begin{equation}
\label{obj_sensr_rng} \tilde{\mathcal{H}}(\mathcal{P}) =
\sum_{i}\int_{(V_{i}\cap\bar{B}(p_i,R))}\phi_n(q)\tilde{f_i}(\|
p_i - q \|))dQ
\end{equation}
where, \[ \tilde{f_i}(r_i) = \begin{cases} f_i(r_i)\quad  \text{if}
\quad r \leq
R_i\\
f_i(R_i)\quad \text{otherwise}
\end{cases}
\]
with $f_i(\cdot) = 1 - \beta_i(\cdot)$ and $\tilde{f}_i(\cdot) = 1 - \tilde{\beta}_i(\cdot)$.

It can be noted that the objective function is made up of sums of
the contributions from sets $V_i \cap \bar{B}(p_i,R_i)$, enabling
the sensors to solve the optimization problem in a spatially
distributed manner.

In reality for range limited sensors the effectiveness should be
zero beyond the range limit. Consider $\hat{f_i}(.) =
\tilde{f_i}(.) - f_i(R_i)$. It can be shown that the objective
function (\ref{obj_sensr_rng}) has the same critical points if
$\tilde{f_i}$ is replaced with $\hat{f_i}$, as the difference in
two objective functions will be a constant term $f_i(R_i)$ (Note
that we have assumed $f_i(R_i) = f_j(R_j) \forall i\text{,}j \in
\{1,\ldots,N\}$.).

\begin{figure}
\centerline{\psfig{figure=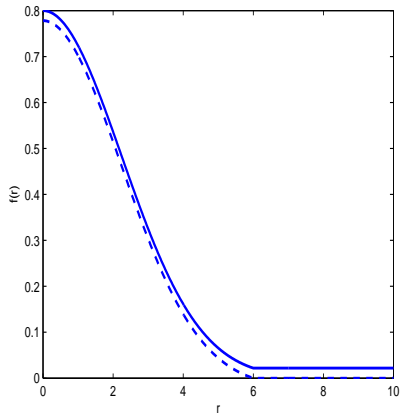,height=6cm,width=6cm}}
\caption{Illustration of $\tilde{f}_i$ and $\hat{f}_i$ in presence of the limit on sensor range. The solid curve represent the sensor effective function $\tilde{f}_i$ and dotted curve is the actual sensor effective function $\hat{f}_i(r) = \tilde{f}_i(r) + (1- f_i(R))$ with $R = 6$.}
\label{sensor_range_fig}
\end{figure}

The gradient of (\ref{obj_sensr_rng}) with respect to $p_i$ can be
determined as
\begin{equation}
\label{grad_sensr_rng}
\frac{\partial(\tilde{\mathcal{H}})}{\partial p_i}(\mathcal{P}) =
2\tilde{M}_{(V_{i}\cap\bar{B}(p_i,R))}(\tilde{C}_{(V_{i}\cap\bar{B}(p_i,R)}
- p_i)
\end{equation}

We use the control law
\begin{equation}
\label{ctrl_sensr_rng} u_i = -k_{prop}(p_{i} -
\tilde{C}_{(V_{i}\cap\bar{B}(p_i,R))})
\end{equation}

\nnd {\it Remark 5.~} It is easy to show, that the gradient
(\ref{grad_sensr_rng}) and the control law (\ref{ctrl_sensr_rng})
are spatially distributed in the $r$-limited (generalized) Delaunay graph
$\mathcal{G}_{LD}$, the Delaunay graph incorporating the sensor
range limits.

\begin{thm}
The trajectories of the sensors governed by the control law
(\ref{ctrl_sensr_rng}), starting from any initial condition
$\mathcal{P}(0)$, will asymptotically converge to the critical
points of $\tilde{\mathcal{H}}$.
\end{thm}

\nnd {\it Proof.~} The proof is similar to that of the Theorem \ref{LaSalle} with $V = - \tilde{\mathcal{H}}(\mathcal{P})$. It can be shown that $V$ is continuously differentiable based on the logic used in Theorem 2.2 in \cite{bullo2}. \hfill $\Box$

Thus, agents having limited sensor range can also be used in HSDS and HCDS. However, as noted earlier, Theorems \ref{DNS_convergence_het} and \ref{CS_convergence_cont_het} are not applicable in this situation as the space spanned by the senors is only a subset of the search space $Q$.

\section{Implementation issues}Here we discuss some of the
theoretical and implementation issues involved in the proposed
search strategies namely, {\em heterogeneous sequential deploy and search} and {\em heterogeneous combined
deploy and search}.

A single step of {\em heterogeneous sequential deploy and search} strategy involves deploying
the agents optimally, and then performing the search task within the
respective Voronoi cells. The deployment step can be
implemented in \emph{continuous time} (as given by control law
(\ref{ctrl1_het})) or in \emph{discrete time} (as in simulations carried
out in this work). When the implementation is in discrete time, in
each time step, the agents move toward the corresponding centroids
and at the end of the deployment step, that is, when the agents are
sufficiently close (as decided by the prescribed tolerance) to the
centroids, the search is performed.

\subsection{Discrete implementation} We convert the
differential equation corresponding to the system dynamics
(\ref{dyn1_het}) to a difference equation.
\begin{equation}
\label{dyna_disc1} \frac{\Delta p_i}{\Delta t} = u_i
\end{equation}
where $\Delta t$ is the discrete time step.

Without loss of generality, we let $\Delta t = 1$ time unit. Then,
(\ref{dyna_disc1}) will be simplified as,
\begin{equation}
\label{dyna_disc} {\Delta p_i}_k = {u_i}_k
\end{equation}
and the control law (\ref{ctrl1_het}) takes the form,
\begin{equation}
\label{ctrl_disc} {u_i}_k = -k_{prop}({p_i}_k - \tilde{C}_{V_{ik}})
\end{equation}
where $k \in \mathbb{N}$ is the iteration count.

\begin{figure*}[ht]
\centerline{
\subfigure[]{\psfig{figure=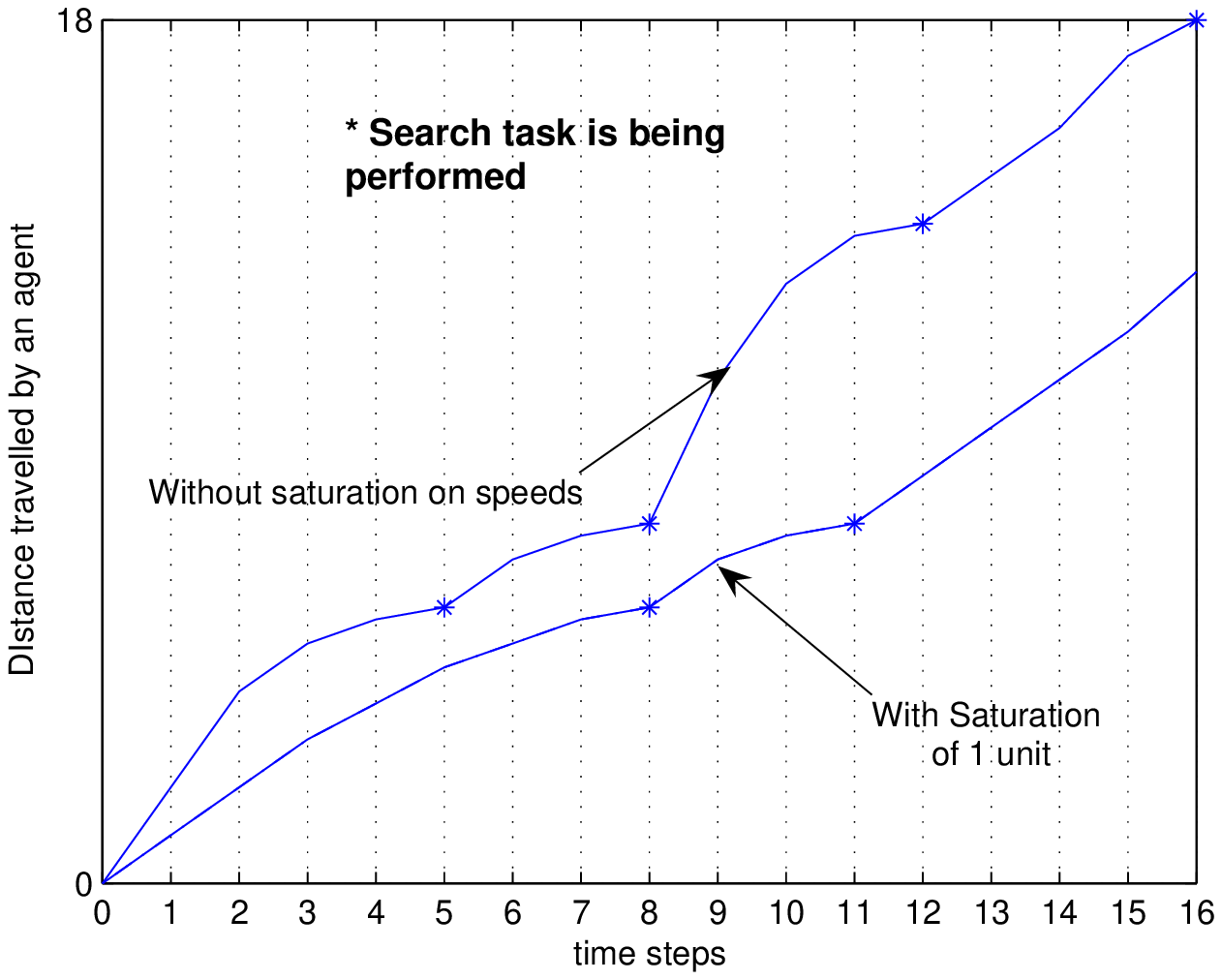,height=6cm,width=6cm}}
\subfigure[]{\psfig{figure=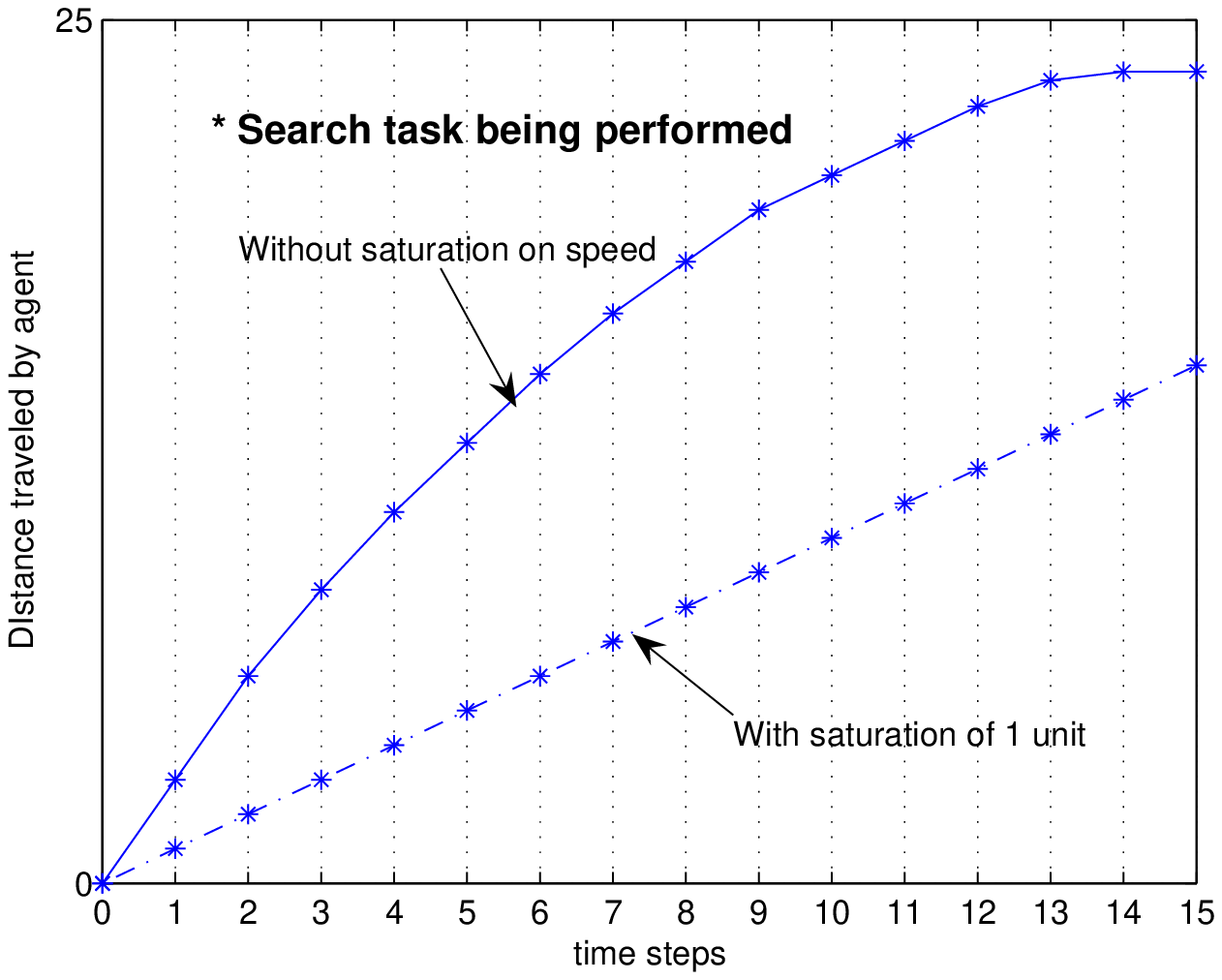,height=6cm,width=6cm}}}
\caption{\small Illustration of the discrete time implementation of
(a) {\em heterogeneous sequential deploy and search} strategy where deployment
takes place according to control law (\ref{ctrl_disc}) at every time
step, and the search task takes place at the end of each deployment
step indicated by `*'and (b) {\em heterogeneous combined deploy and search}
strategies where search is performed at every step while the agents
move according to the control law (\ref{ctrl_disc}).}
\label{dns_step_illustrate}
\end{figure*}

The control input ${u_i}_k$ is the desired speed of the $i$-th agent
at the $k$-th time step. The agent moves with this speed for $\Delta
t$ time units. With $\Delta t = 1$,  ${u_i}_k$ acts as increment of
$p_i$ per step. In other words ${u_i}_k = \Delta p_i = {p_i}_{k+1} -
{p_i}_k$. Thus, if $\Delta t = T$ time units, then the search task
takes place after $mT$ time units, where $m$ is a non-zero integer,
the number of time steps taken to achieve the optimal deployment.
The process is illustrated in Figure \ref{dns_step_illustrate}(a). Search is performed at the end of deployment step and the time instances at which search is performed are marked with `*'. The consecutive search task is performed at a time interval of at
most $T$ time units. We define the {\em latency}, $t_s$, of the
agents as the maximum time taken to acquire the information, process
it, and successfully update the uncertainty density. Finally, $T$
should be chosen to be greater than or equal to $t_s$.

In the {\em heterogeneous combined deploy and search} strategy the agents perform
search while moving toward corresponding centroids, without
waiting till the end of the deployment step. The agents perform search after every $T$ time units
as they move toward the respective centroids according to the
discrete control law (\ref{ctrl_disc}). The process is illustrated
in Figure \ref{dns_step_illustrate}(b). Search is performed at every time instance as indicated by `*'. The plot of distance traveled versus time steps is smoother in case of CDS compared to that of SDS. 

\subsection{Effect of saturation}
The control inputs given by control law (\ref{ctrl1_het}) or
(\ref{ctrl_disc}) determine the speeds of the agents. In practical
implementation, it is likely that there will be a constraint on the
maximum speed of agents. Such a limit will appear as a saturation on
the control input. In case of the {\em heterogeneous sequential deploy and search} strategy,
the effect of saturation on control input might lead to slower
convergence of the deployment step. During the initial few steps, it
is likely that the control input provided by (\ref{ctrl1_het}) can cross
the saturation limit, whereas later, as the agents approach the
centroids, the control input naturally reduces (as it is
proportional to the distance between the agents and the respective
centroids). Thus, effect of saturation is at most a possible
increase in the time gap between consecutive search steps as
illustrated in Figure \ref{dns_step_illustrate}(a).

In the case of {\em heterogeneous combined deploy and search} strategy,
whenever the control effort computed by (\ref{ctrl1_het}) crosses the
saturation limit, the actual speed is limited to the saturation
value. The time lag between two consecutive search tasks remains
fixed at $T$ irrespective of the saturation. But, with saturation,
the distance traveled by the agents before performing the next
search task reduces (or remains the same if the control input given
by (\ref{ctrl1_het}) is less than the saturation limit). This is
illustrated in Figure \ref{dns_step_illustrate}(b). This will
probably result in a faster search due to frequent searches.

\subsection{Spatial distributedness}
Here we discuss the implication of spatial distributedness of the proposed search strategies from a practical
point of view. We have seen that both the search strategies are spatially distributed in the generalized Delaunay graph. These
results imply that all the agents need to do computations based on only local information, that is, by the knowledge
about position of neighboring agents. Also, the agents should have access to the updated uncertainty map within
their Voronoi cells. This can be achieved in several ways. One such way is that all the agents communicate with
a central information provider. But it is not necessary to have this global information. Each agent can communicate with its Voronoi neighbors ($\mathcal{N}_{\mathcal{G}}(i)$) and obtain the updated uncertainty information in a region $\cup_{\mathcal{N}_{\mathcal{G}}(i)} V_i$. As the generalized  Voronoi partition $\{V_i\}$ depends at least continuously on $\mathcal{P}$, the agent configuration, in an evolving generalized Delaunay graph, the communication within the neighbors is sufficient for each agent to obtain the uncertainty within its new Voronoi cell. The issues related
to communication of uncertainty information are not addressed in the paper except to assume that uncertainty
information is available to the agents. It is also possible that the agents can estimate the uncertainty map as done
in \cite{schwager}.

In practical conditions, the agents can communicate with other agents only when they are within the limits of
the sensor range. The generalized Delaunay graph does not allow sensor range limitations to be incorporated. We need to use
\emph{$r$-limited generalized Delaunay} graph or \emph{$r$-generalized Delaunay} graph (generalized versions of \emph{$r$-limited Delaunay} graph and \emph{$r$-Delaunay} \cite{graph}) to incorporate the sensor range limitations. It needs to be studied if
the proposed search strategies are still spatially distributed on these graphs. In any case, the scenario changes with
incorporation of sensor range limitations into the search strategies. The updating of uncertainty density will also be
within the sensor range limits (in fact, it is within the region $V_i \cap \bar{B}(p_i,R)$). The centroid that is computed will also be within the new restricted area. For an optimal
deployment problem, from the perspective of sensor coverage, it has been observed that the corresponding control law
is still spatially distributed (in $r$-limited generalized Delaunay graph) and globally asymptotically stable.

\subsection{Synchronization}
Synchronization plays an important role in multi-agent systems. Here we discuss this issue for both HSDS and
HCDS strategies. In the case of HSDS, theoretically all the agents reach the respective centroid at infinite time. But in a
practical implementation, the agents are required to be sufficiently close, where the closeness is suitably defined, to
the respective centroids before starting the search operation. It is possible that at any point in time, different agents
are at different distances from the corresponding centroid. The agents need to come to a consensus as to when to
end the deployment and perform search operation. We have implemented the strategy in a single centralized program using MATLAB. In a practical situation, synchronization can be attained by agents communicating a flag bit indicating if an agent has reached its centroid or not. When all the agents  have reached the respective centroid within a tolerance distance, the search can be performed. We also assume that sensing and communication are instantaneous. In our simulation experiments we assume such a communication exists. Since the objective of this work is to evaluate the effectiveness of the search strategies, we make assumptions that simplify implementation without affecting the search effectiveness.

HCDS operates in a synchronous manner by design. If all the agents start at the same instant of time and have
synchronized clocks, the search task is performed by every agent after the same interval of time. Given an accurate
global clock, synchronization is not a major issue in case of HCDS.

Further in [17] authors provide an asynchronous implementation for coverage control which can be suitably modified for HSDS and HCDS.

\section{Results and Discussions}
A
few simulation experiments were carried out to validate the
proposed heterogeneous multi-agent search strategies. We used
$\beta_i = 1 -k_ie^{-\alpha_ir_{i}^2}$, $i \in \{1, \ldots, N\}$.
The search space $Q$ is a square area of size $10 \times 10$ in
$\mathbb{R}^2$. We vary the parameters $k_i$s and $\alpha_i$s to
get different cases.

\begin{figure}
\centerline{
\subfigure[]{\psfig{figure=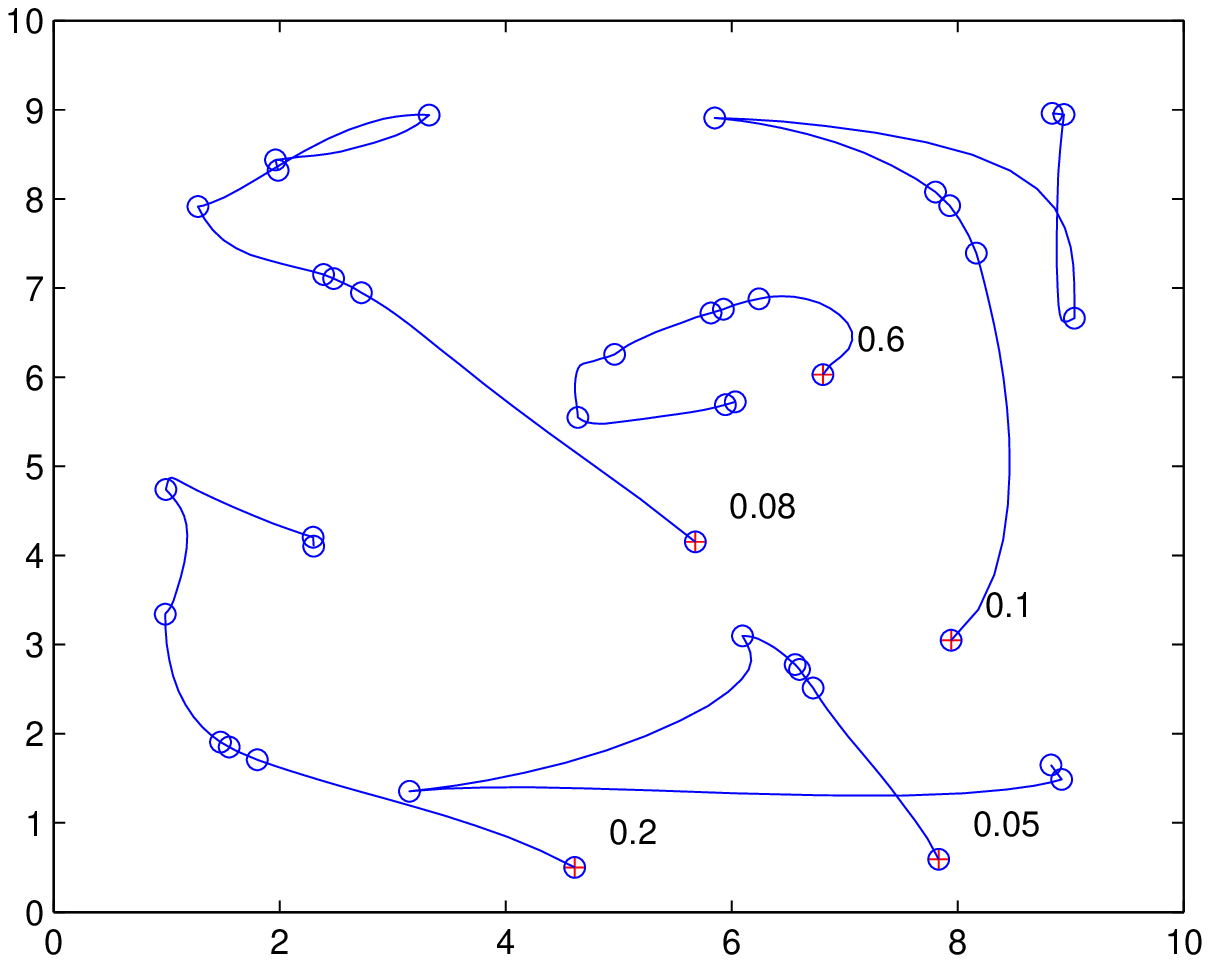,width=6cm,height=6cm}}
\subfigure[]{\psfig{figure=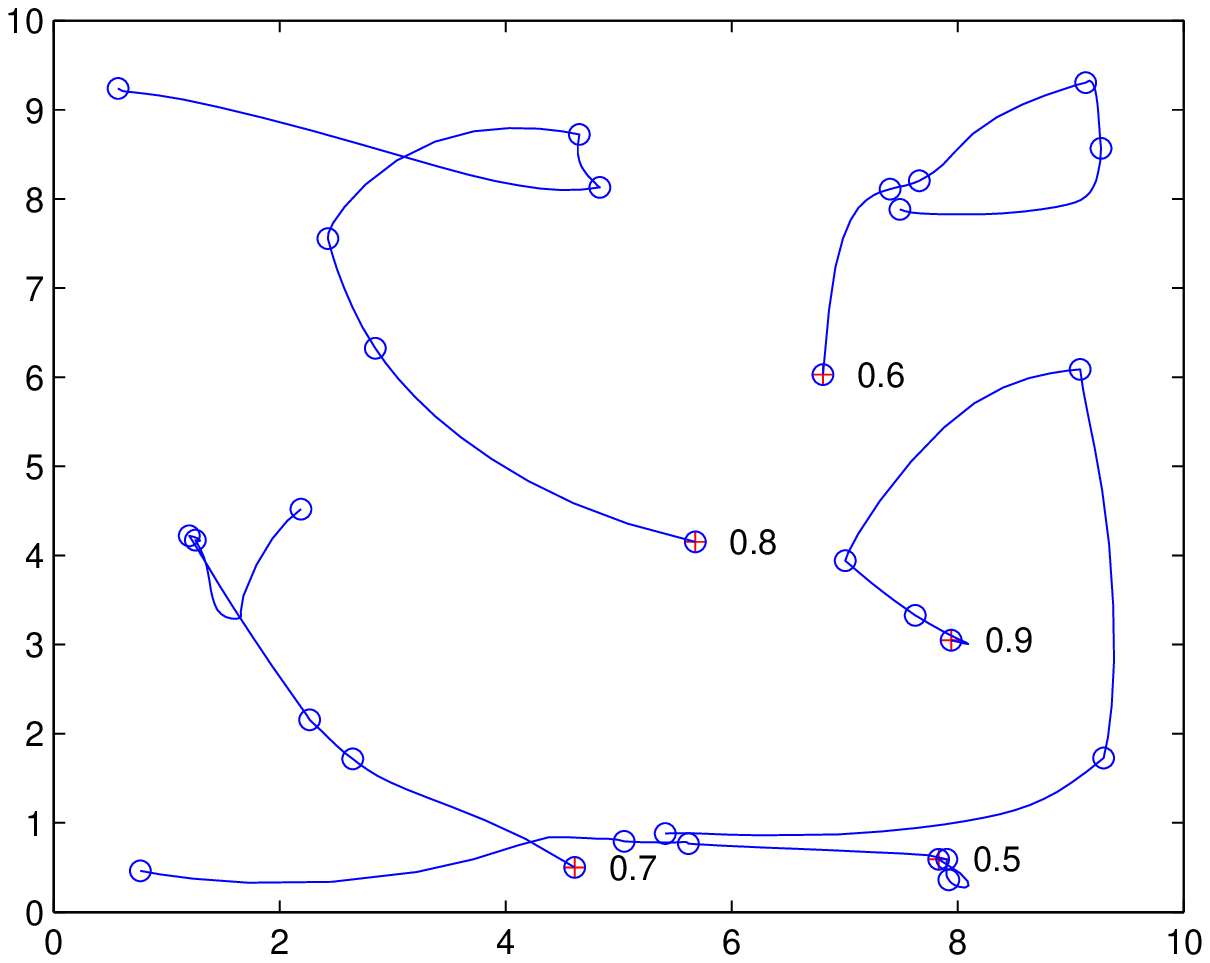,width=6cm,height=6cm}}
\subfigure[]{\psfig{figure=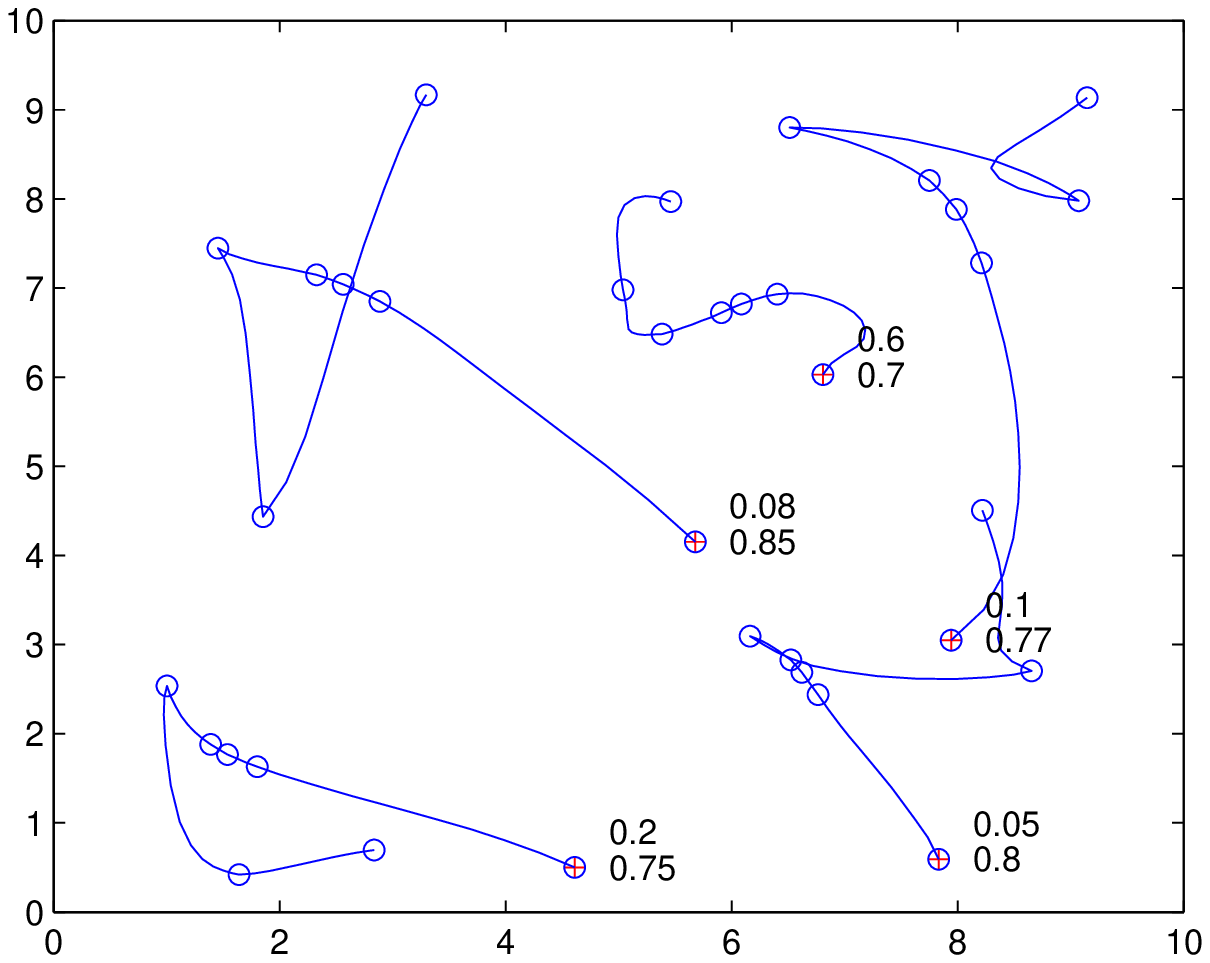,width=6cm,height=6cm}}}
\centerline{
\subfigure[]{\psfig{figure=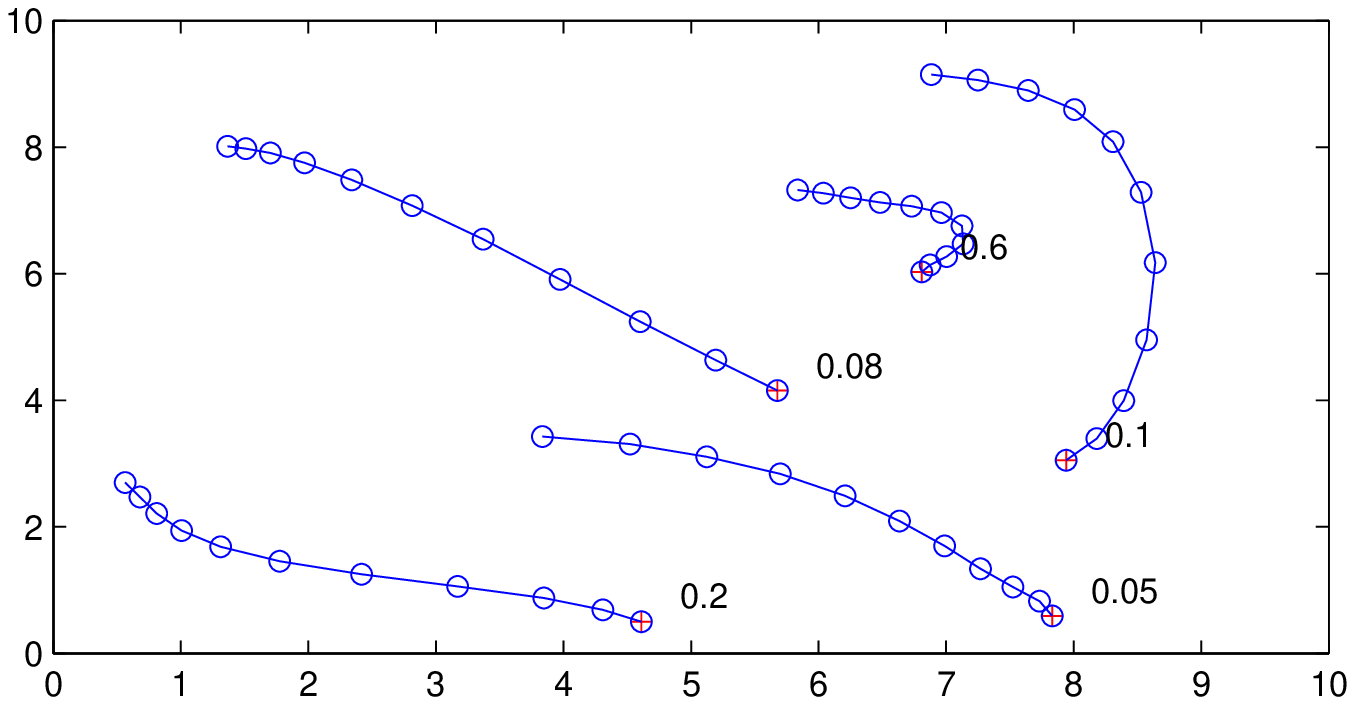,width=6cm,height=6cm}}
\subfigure[]{\psfig{figure=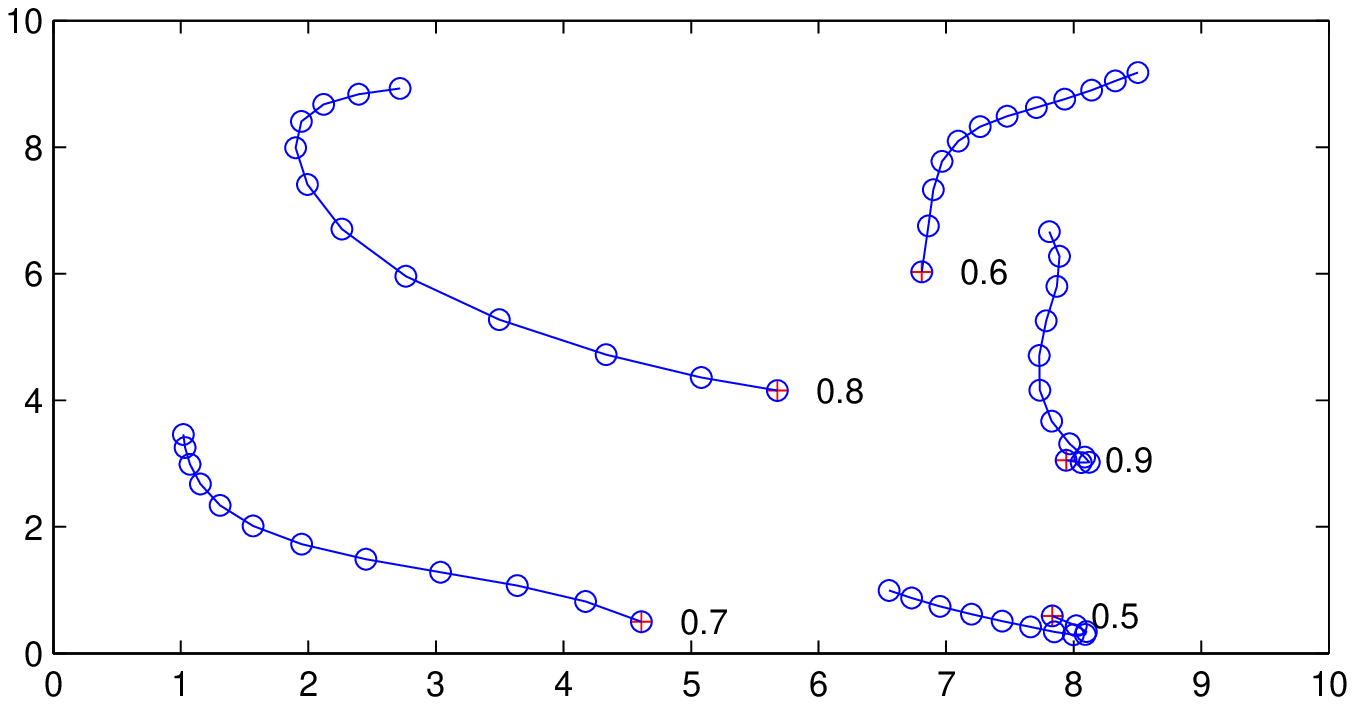,width=6cm,height=6cm}}
\subfigure[]{\psfig{figure=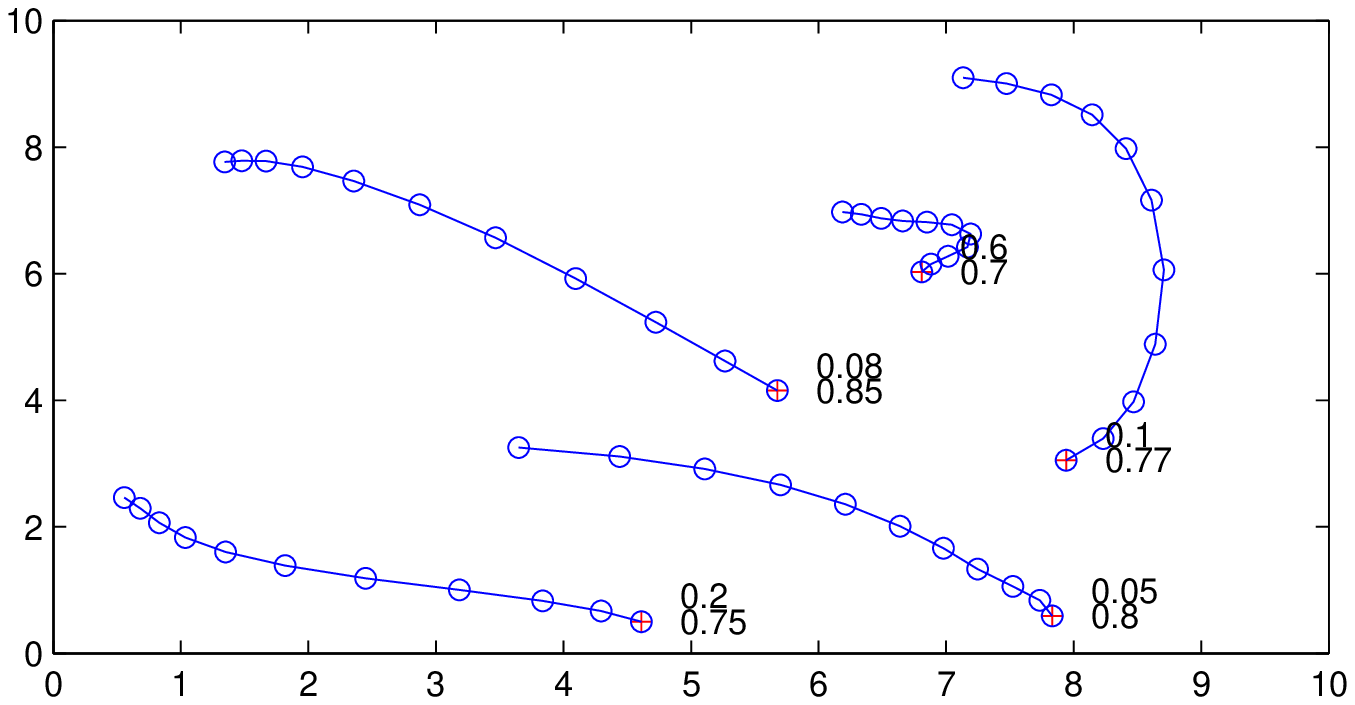,width=6cm,height=6cm}}}
\caption{Trajectories for agents with HSDS for $\beta_i(r_i) = 1-k_i e^{\alpha_i r^2_i}$ with (a) $\alpha_i$ varying, (b) $k_i$ varying, and (c) both $\alpha_i$ and $k_i$ varying. Corresponding trajectories with HCDS are shown in (d), (e), and (f). `o' indicate the search
instances. The parameters are indicated near the starting location of
corresponding agents, which are indicated by `+'.}\label{trajfig}
\end{figure}

\begin{figure}
\centerline{
\subfigure[]{\psfig{figure=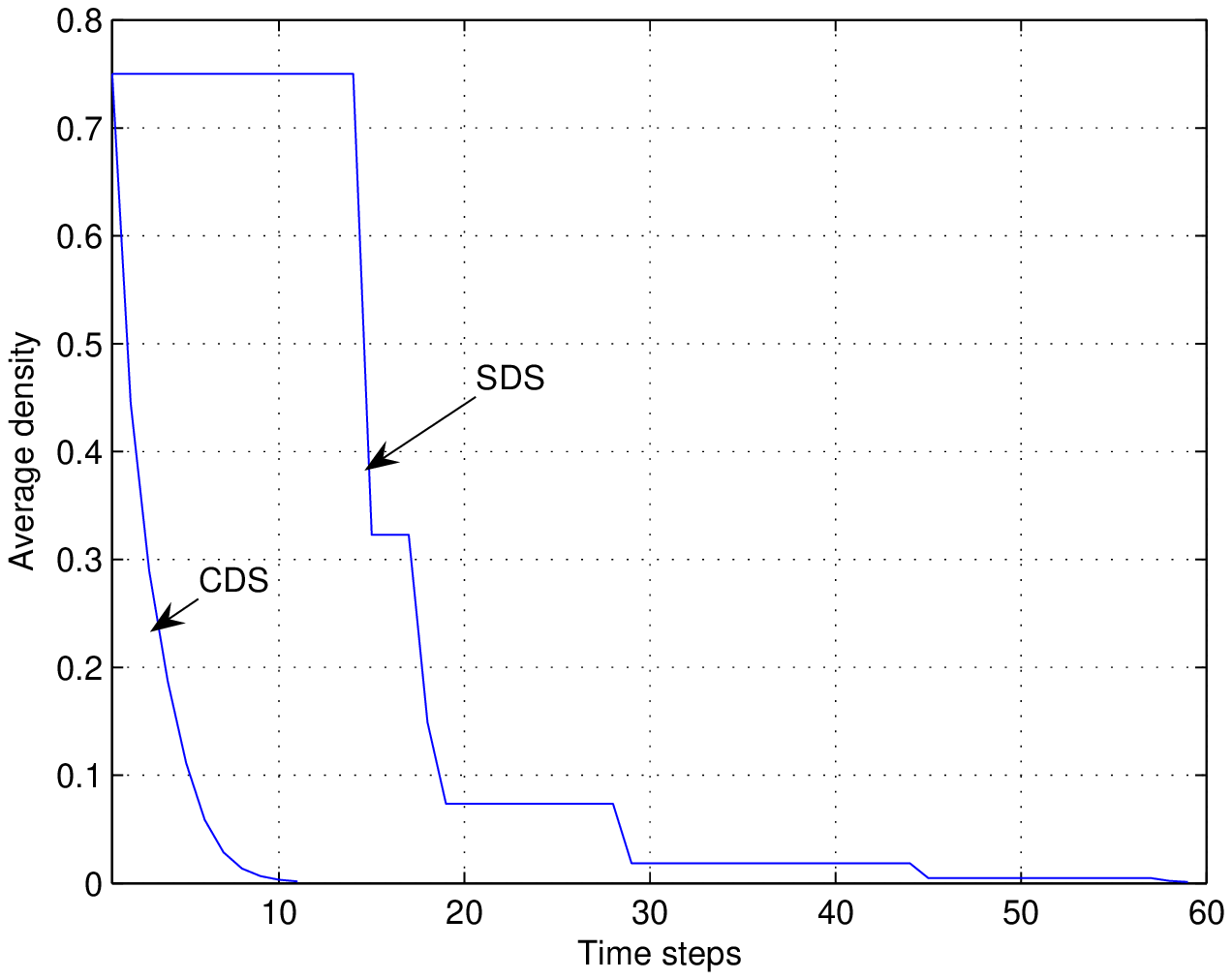,width=6cm,height=6cm}}
\subfigure[]{\psfig{figure=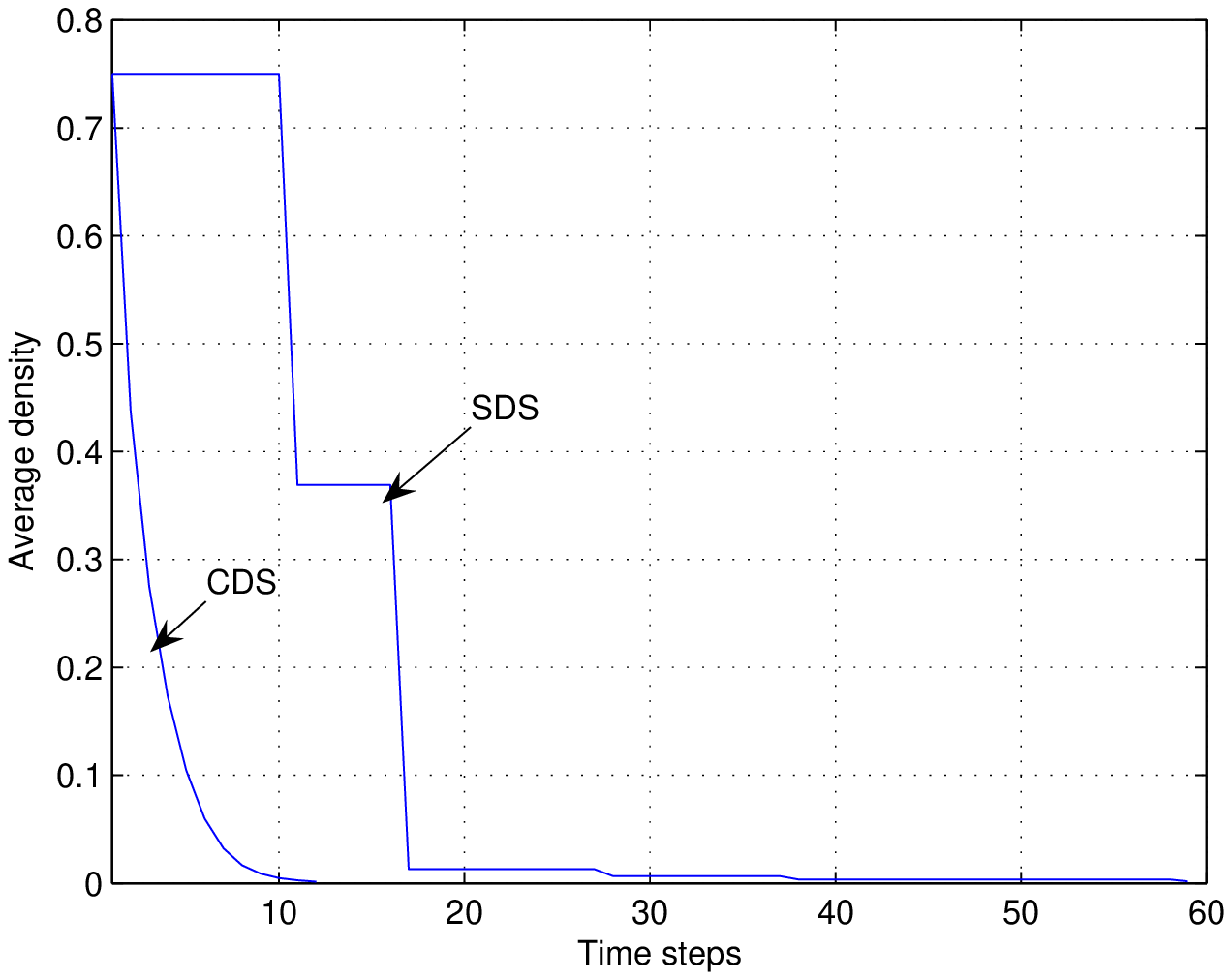,width=6cm,height=6cm}}
\subfigure[]{\psfig{figure=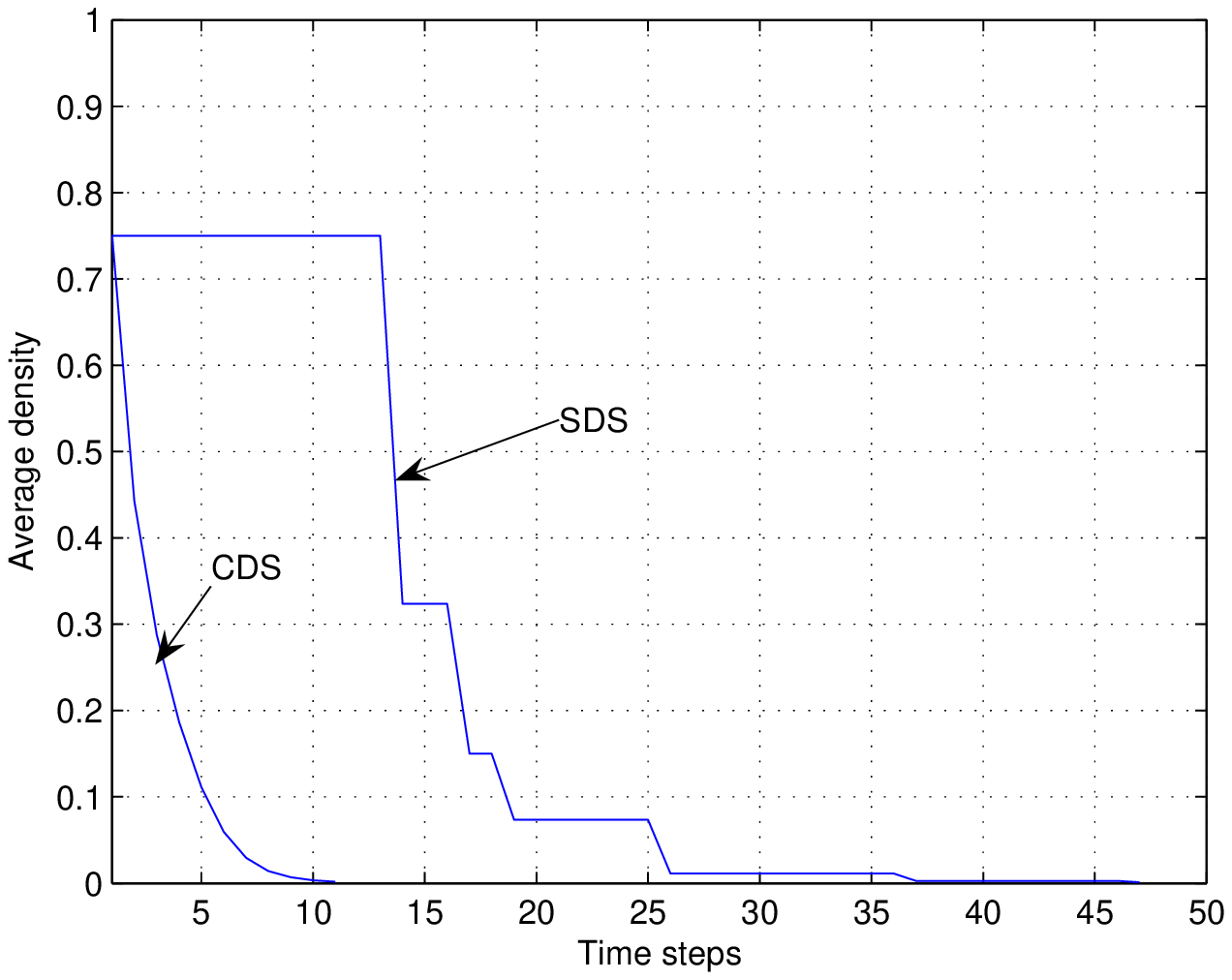,width=6cm,height=6cm}}
}
\centerline{
\subfigure[]{\psfig{figure=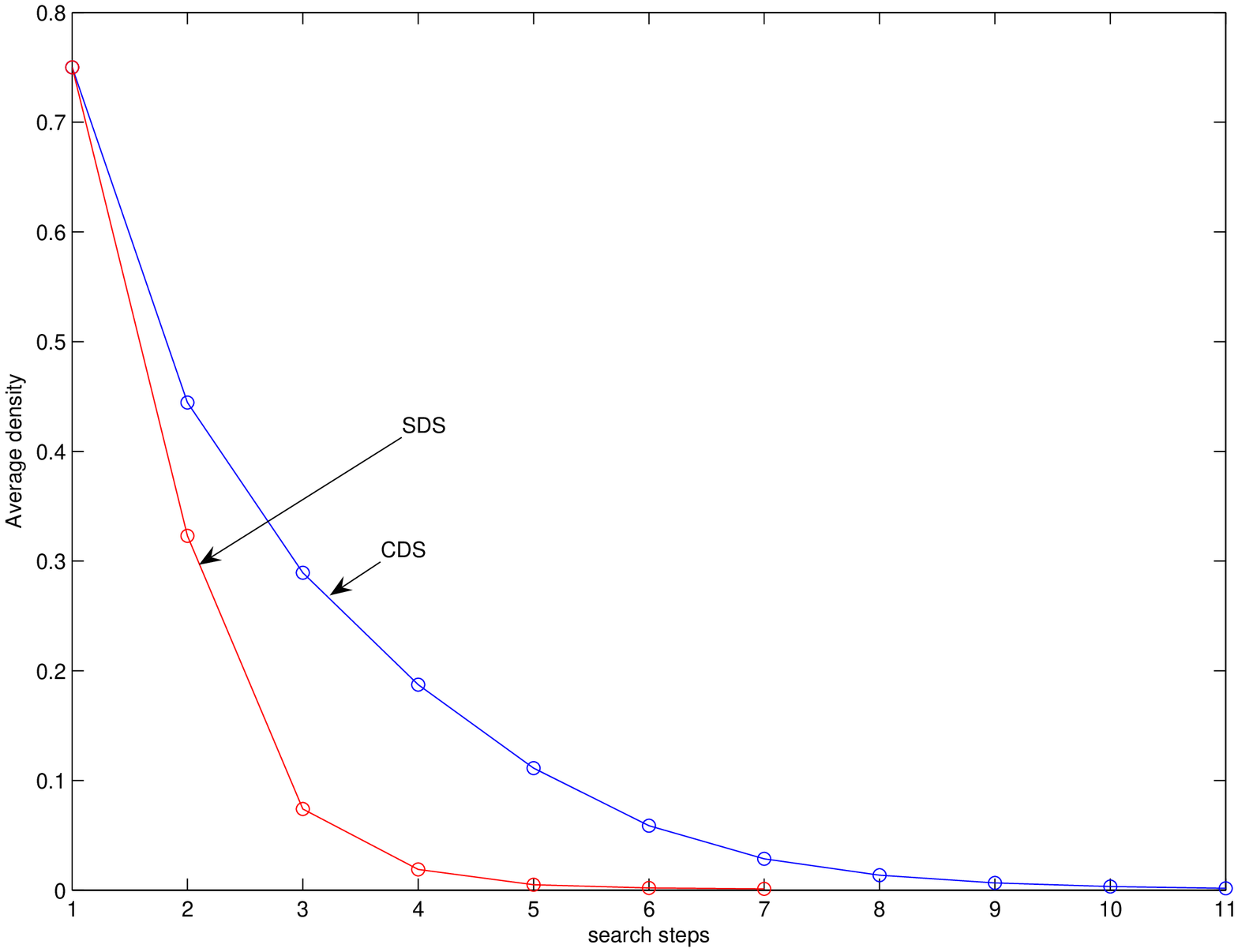,width=6cm,height=6cm}}
\subfigure[]{\psfig{figure=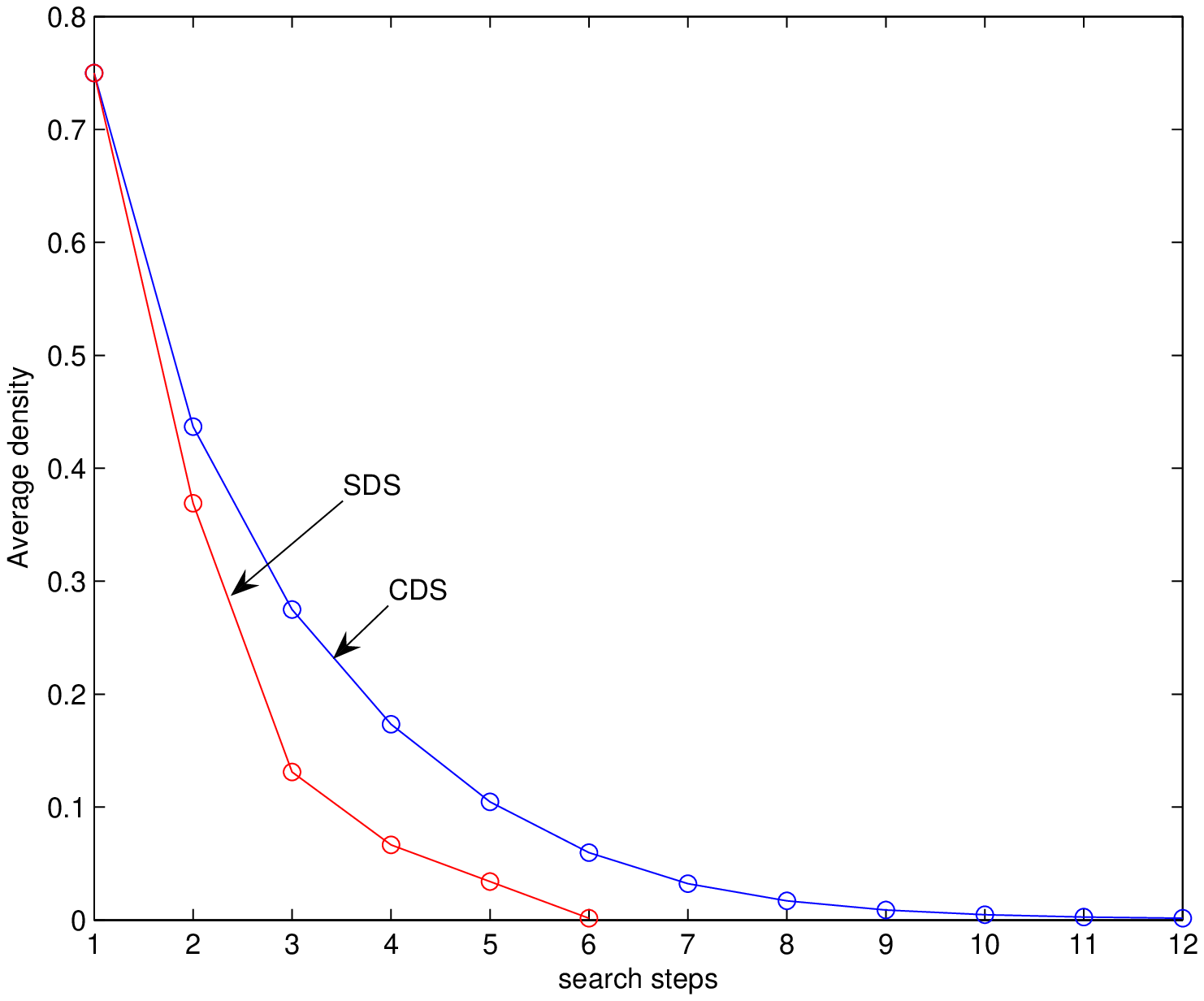,width=6cm,height=6cm}}
\subfigure[]{\psfig{figure=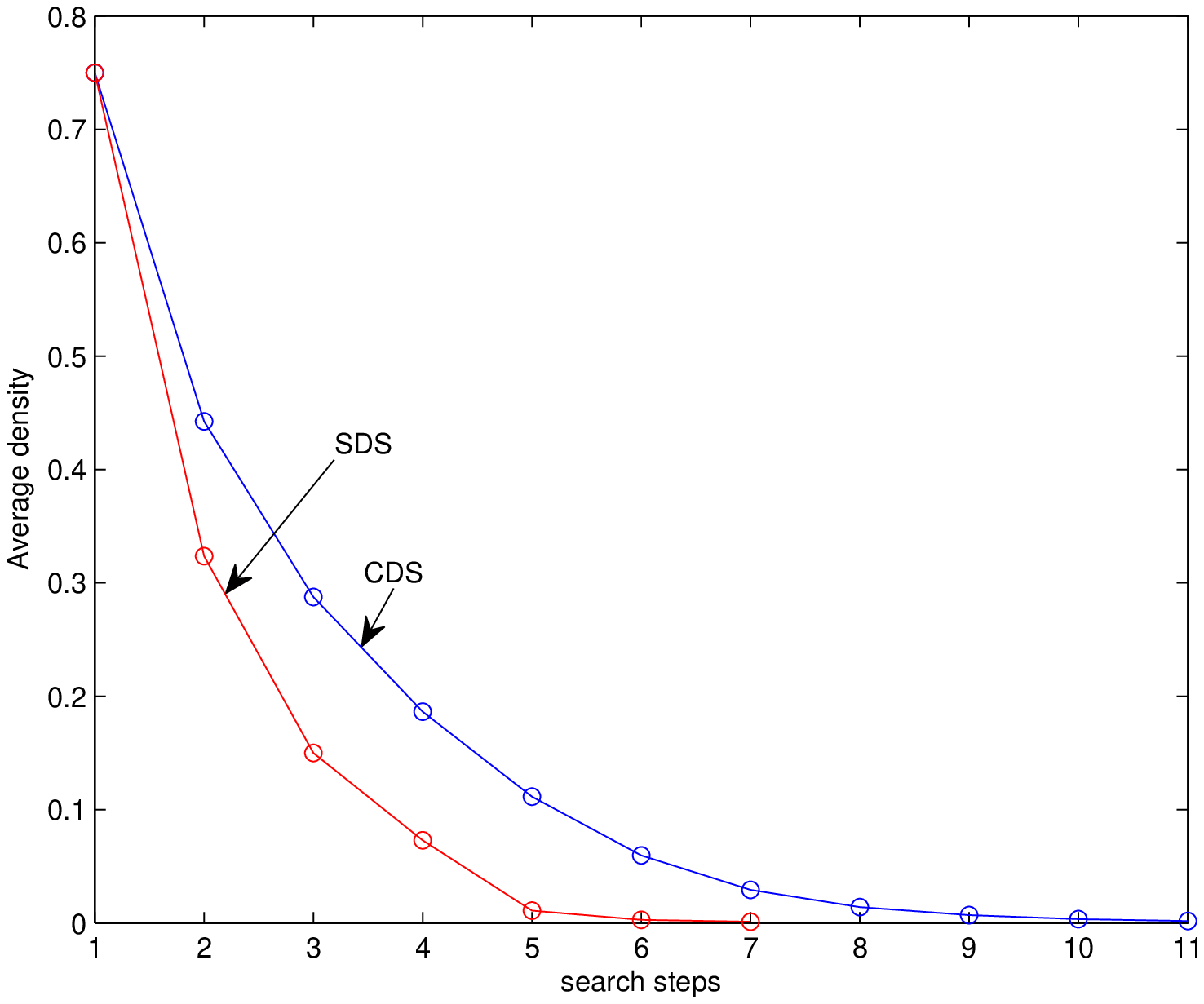,width=6cm,height=6cm}}
 }
\caption{Reduction of uncertainty with HSDS and HCDS with time steps for $\beta_i(r_i) = 1-k_i e^{\alpha_i r^2_i}$ with (a) $\alpha_i$ varying,
(b) $k_i$ varying, and (c) both $\alpha_i$ and $k_i$ varying. Corresponding reduction in uncertainty with number of search steps are shown in (d), (e),
and (f). `o' indicate the search instances.} \label{densfig}
\end{figure}

Figures \ref{trajfig} and \ref{densfig} show results of simulation
experiments carried out. Figure \ref{trajfig}(a) and (d) show the
trajectories of the agents for HSDS and HCDS strategies, while
Figure \ref{densfig}(a) shows the average uncertainty density
history for two strategies with $k_i=0.8$ and different
$\alpha_i$. Figure \ref{densfig} (d) shows the average uncertainty density with number of searches. The starting location of the agents are marked with `+' and end of each of the deploy steps are marked with `o' along the
trajectories for the case of HSDS (Figure \ref{trajfig}(a)). It can be observed from Figure \ref{densfig} (a), that HCDS reduces uncertainty much faster in time compared to HSDS. HSDS takes requires about 20 time steps to reduce the uncertainty below a value of 0.8, while HCDS achieves same reduction in about 6 time steps. This is due to increased frequency of searches in case of HCDS compared to that of HSDS. But Figure \ref{densfig} (d) reveals that HSDS performs better than HCDS in terms of requiring lesser number of search steps for same amount of uncertainty reduction. In case of HSDS, the agents get optimally deployed before performing the search and hence, compared to HCDS, the uncertainty density is higher in each search step of HSDS.

Figures \ref{trajfig} (b) and (e) show trajectories of agents with $\alpha_i = 0.1$ and varying $k_i$ for HSDS and HCDS, and \ref{densfig} (b) and (e) show the average uncertainty with number of time steps and number of search steps for both the strategies. Figures \ref{trajfig} (c) and (f), and \ref{densfig} (c) and (f) show corresponding results with both the parameters $k_i$ and $\alpha_i$
varying. Values of the varying parameters have been indicated in
the figures, near the starting point of each agent which marked with `+'.

In all the cases, it can be observed that the trajectories of
agents with HCDS strategy are considerably smoother and shorter
than those corresponding to the HSDS strategy, nevertheless, both
strategies successfully reduce the uncertainty density. The
simulation results also indicate that the HCDS strategy performs
better even in terms of faster reduction of the average
uncertainty density, while HSDS  performs better in terms of requiring fewer search instances. It can also be observed that the agents move away from each other covering the search space in a cooperative manner. Figures \ref{densfig} (a), (b), and (c) illustrate that in HCDS, the search is performed at every time instance, and in HSDS, the search is performed only after optimal deployment of agents.

The simulation results demonstrate that both the proposed
heterogeneous search strategies perform well as indicated by the
theoretical analysis and that the HCDS strategies performs well in
terms of shorter and smoother agent trajectories and faster uncertainty reduction.

When the sensors have heterogeneous capabilities either in terms
of effectiveness (indicated by the parameter $k_i$) or in terms of
the range (indicated by the parameter $\alpha_i$), those with
higher effectiveness share more load, while the weaker ones may
remain inactive (in extreme cases of heterogeneity). 
The motivation here is to demonstrate the ability to handle
heterogeneity in sensors' capabilities. The sensor parameters are
assumed to be given and the problem is of designing a suitable
search strategy. It might be interesting to select sensors with
required parameters so as to improve the effectiveness of the search
strategy. But this design optimization is beyond the scope of this
work.

\section{Conclusions} We have used a generalization of the
Voronoi partition  to formulate and solve a heterogeneous
multi-agent search problem. The agents having sensors with
heterogeneous capabilities were deployed in the search space in an
optimal way maximizing per step search effectiveness. The objective
function, its critical points, a control law that determines the
agent trajectory, its spatial distributedness and convergence
properties were discussed. Based on the optimal deployment strategy,
two heterogeneous multi-agent search strategies namely {\em heterogeneous
sequential deploy and search} and {\em heterogeneous combined deploy and search}
have been proposed and their spatial distributedeness and
convergence properties have been studied. Effect of constraints on agents' speeds and limit on sensor range have been discussed.
The simulation experiments demonstrate that the search strategies perform well. The {\em heterogeneous combined
deploy and search} strategy is seen to perform better in terms of
shorter, smoother agent trajectory and faster search.

Analysis of the properties of the generalized Voronoi partition is one of the possible direction for research. Work on effective algorithms for computation related to generalized Voronoi partition will be very useful in effective implementation of the search strategies presented in this paper.  Further generalization of the Voronoi partition so as to incorporate anisotropy in the sensors along with the heterogeneity can also be a very useful exercise. With such a generalization, search strategies for multiple agents equipped with heterogeneous and anisotropic sensors can be formulated.
It is also interesting to explore new applications of the generalized Voronoi partition.

\end{document}